\documentclass[11pt]{amsart}
\usepackage{a4,amsmath,amssymb,amscd,amsfonts,stmaryrd,verbatim}

\newif\ifnige
\nigefalse

\ifnige

  \def\url#1{\verb=#1=}
  \def\arxiv#1{\texttt{arXiv:#1}}

\else
  \usepackage{lmodern}
  \usepackage{stmaryrd}
  \usepackage{mathtools}
  \usepackage[pdfusetitle]{hyperref}

  \def\arxiv#1{\href{http://arxiv.org/abs/#1}{\texttt{arXiv:#1}}}

\fi
  
\newcommand{\lebn}{\tt 19.08.11 \jobname}

\newtheorem{prop}[equation]{Proposition}

\newtheorem{thm}[equation]{Theorem}
\newtheorem{cor}[equation]{Corollary}
\newtheorem{lem}[equation]{Lemma}

\theoremstyle{definition}

\newtheorem{defn}[equation]{Definition}
\newtheorem{defns}[equation]{Definitions}
\newtheorem{rem}[equation]{Remark}
\newtheorem{rems}[equation]{Remarks}

\newtheorem{exa}[equation]{Example}

\theoremstyle{remark}
\newtheorem*{acknowledgements}{Acknowledgements}

\numberwithin{equation}{section}

\newcommand{\spandsp}{\mbox{$\qquad\text{and}\qquad$}}

\newcommand{\sands}{\mbox{$\quad\text{and}\quad$}}

\newcommand{\sts}[1]{\mbox{$\quad\text{#1}\quad$}}

\newcommand{\Hom}{\operatorname{Hom}}

\ifnige
  \newcommand{\letbe}{\mathbin{\raisebox{.3pt}{:}\!=}}
\else
  \def\letbe{\vcentcolon =}
\fi
\newcommand{\br}[1]{\mbox{$\langle #1\rangle$}}
\newcommand{\bbr}[1]{\mbox{$\big\langle #1\big\rangle$}}

\newcommand{\?}{\mskip-2mu}

\newcommand{\3}{\mskip1mu}
\newcommand{\pfm}{\mskip.5mu}
\newcommand{\scirc}{\mathbin{\raisebox{.6pt}{$\scriptstyle \circ$}}}

\newcommand{\CPn}{\mbox{$\bC P^n$}}

\newcommand{\ssp}{{\scriptscriptstyle{+}}}

\newcommand{\MU}{\mbox{\it MU}}

\newcommand{\osu}{\mbox{$\varOmega_*^{U}$}}
\newcommand{\ous}{\mbox{$\varOmega^*_{U}$}}

\newcommand{\KO}{\mbox{\it KO\/}}

\newcommand{\bC}{\mathbb{C}}

\newcommand{\bN}{\mathbb{N}}
\newcommand{\bP}{\mathbb{P}}
\newcommand{\bQ}{\mathbb{Q}}
\newcommand{\bR}{\mathbb{R}}
\newcommand{\bZ}{\mathbb{Z}}

\newcommand{\srf}[1]{\mbox{${\text{\it SR}^F}$}}
\newcommand{\Th}{\mbox{\it Th}}

\newcommand{\daaja}{Davis and Januszkiewicz}

\newcommand{\Pc}{\mbox{$\bP(\chi)$}}
\newcommand{\Pd}{\mbox{$\bP(\chi')$}}
\newcommand{\Ppc}{\mbox{$\bP({}_p\chi)$}}
\newcommand{\Pp}{\mbox{$\bP(\pi)$}}

\newcommand{\Po}{\mbox{$\bP(\omega)$}}

\newcommand{\cat}[1]{\mbox{\sc #1}}

\newcommand{\ra}{\rightarrow}
\newcommand{\lra}{\longrightarrow}
\newcommand{\llra}{\relbar\joinrel\hspace{-2pt}\lra}
\newcommand{\lllra}{\relbar\joinrel\hspace{-2pt}\llra}
\newcommand{\llllra}{\relbar\joinrel\hspace{-2pt}\lllra}
\newcommand{\lllllra}{\relbar\joinrel\hspace{-2pt}\llllra}
\newcommand{\la}{\leftarrow}
\newcommand{\lla}{\longleftarrow}
\newcommand{\llla}{\lla\joinrel\hspace{-2pt}\relbar}
\newcommand{\lllla}{\llla\joinrel\hspace{-2pt}\relbar}
\newcommand{\llllla}{\lllla\joinrel\hspace{-2pt}\relbar}

\makeatletter \newcommand{\Ddots}{\mathinner{\mkern1mu\raise\p@
\vbox{\kern7\p@\hbox{.}}\mkern2mu
\raise4\p@\hbox{.}\mkern2mu\raise7\p@\hbox{.}\mkern1mu}}
\makeatother

\newcommand{\lcm}{\operatorname{lcm}}

\begin{document}
\bibliographystyle{plain}

\title[Weighted projective spaces and iterated Thom spaces]
{Weighted projective spaces and\\ iterated Thom spaces}

\author{Anthony Bahri}
\address{Department of Mathematics, Rider University, Lawrenceville,
NJ 08648, USA}
\email{bahri@rider.edu}

\author{Matthias Franz}
\address{Department of Mathematics, University of Western Ontario,
London, ON N6A~5B7, Canada}\email{mfranz@uwo.ca}

\author{Nigel Ray}
\address{School of Mathematics, University of Manchester, Oxford
Road, Manchester M13~9PL, England}
\email{nigel.ray@manchester.ac.uk}

\ifnige\else
\pdfstringdefDisableCommands{\def\and{, }}
\hypersetup{pdfauthor=\authors}
\fi

\thanks{A\,B was partially supported by a Rider University Summer 
Research Fellowship and Grant~\#210386 from the Simons Foundation; 
M\,F was partially supported by an NSERC Discovery Grant; N\,R was 
partially supported by Royal Society International Joint Projects 
Grant \#R104514.}
\keywords {Complex oriented cohomology theory, iterated Thom space,
Thom isomorphism, weighted projective space} 

\subjclass[2010]{Primary 57R18; secondary 14M25, 55N22}

\begin{abstract}

  For any $(n+1)$--dimensional weight vector $\chi$ of positive
  integers, the weighted projective space $\bP(\chi)$ is a projective
  toric variety, and has orbifold singularities in every case other
  than $\bC P^n$. We study the algebraic topology of $\bP(\chi)$,
  paying particular attention to its localisation at individual primes
  $p$. We identify certain $p$-primary weight vectors $\pi$ for which
  $\bP(\pi)$ is homeomorphic to an iterated Thom space over $S^2$, and
  discuss how any $\bP(\chi)$ may be reconstructed from its
  $p$-primary factors. We express Kawasaki's computations of the
  integral cohomology ring $H^*(\bP(\chi);\bZ)$ in terms of iterated
  Thom isomorphisms, and recover Al Amrani's extension to complex
  $K$-theory. Our methods generalise to arbitrary complex oriented
  cohomology algebras $E^*(\bP(\chi))$ and their dual homology
  coalgebras $E_*(\bP(\chi))$, as we demonstrate for complex cobordism
  theory (the universal example). In particular, we describe a
  fundamental class in $\varOmega_{2n}^U(\bP(\chi))$, which may be
  interpreted as a resolution of singularities.

\end{abstract}

\maketitle

%
%
%
%
%
%
%
%
%

\section{Introduction}\label{intro}

Weighted projective spaces $\Pc$ are defined for every integral
weight vector $\chi$, and constitute a family of singular toric
varieties on which many hypotheses may be tested. Our central interest
is to study their algebraic topology, and identify them up to
homeomorphism with spaces whose singularities are more familiar,
namely iterated Thom complexes. In the process of working towards this
goal, it transpires that the cohomology algebras $E^*(\Pc)$, and
homology coalgebras $E_*(\Pc)$, may be described in terms of iterated
Thom isomorphisms for any complex oriented cohomology and homology
theories $E^*(-)$ and $E_*(-)$. This provides a fruitful new
perspective on algebraic objects of considerable complexity, and we 
present several explicit examples below.

Until the last few years, literature on the algebraic topology of
weighted projective spaces has been sparse, and restricted mainly to
work of Kawasaki \cite{kaw:ctp} and Al Amrani \cite{ala:csw},
\cite{ala:cbc}.  Immediately after the latter, Nishimura--Yosimura
\cite{niyo:qkt} took up the challenge of computing the real $K$-theory
groups $\KO^*(\Pc)$, whose difficulty is increased by the lack of
complex orientation for $\KO^*(-)$. More recently, it has become
apparent that the toric structure of $\Pc$ is particularly important,
and our own work \cite{bafrra:ecr}, \cite{bafrnora:cwp} has exploited
this fact. 
Most of our current results are independent of the toric framework,
but its beauty is sufficiently compelling that we have retained here
it for background and motivation.

We present our results in the following order.

Section \ref{wps} establishes notation by recalling definitions of
$\Pc$ and the associated weighted lens spaces $L(k;\chi)$ in the toric
setting. A crucial cofibre sequence of Kawasaki is expressed in terms
of normal neighbourhoods over the orbit simplex $\varDelta^n$ of the
canonical torus action, and his resulting computation of the integral
cohomology rings $H^*(\Pc)$ and $H^*(L(k;\chi))$ are summarised. The
important maps $\phi\colon\CPn\ra\Pc$ and $\psi\colon\Pc\ra\CPn$ are
described in terms of homogeneous coordinates, and the complex line
bundle classified by $\psi$ is identified for future use.

Section \ref{itthspprpa} introduces the concept of \emph{divisive}
weight vectors and their normalisations, and Theorem \ref{divisprop}
employs an observation of Al Amrani (relating $\psi$ and normal
neighbourhoods) to express $\Pc$ as an iterated Thom complex whenever
$\chi$ is divisive. Fundamental examples are provided by
\emph{$p$-primary} weight vectors $\pi$, whose coordinates are powers
of a fixed prime $p$. Archetypal examples are the \emph{$p$-contents}
${}_p\chi$ of $\chi$, which 
lead to the $p$-primary decomposition of $\chi$. Maps
$e(\chi/\omega)\colon\Po\ra\Pc$ are introduced as common
generalisations of $\psi$ and $\sigma$, and determine
\emph{extraction} and \emph{inclusion} maps between $\Pc$ and its
$p$-primary parts $\Ppc$. The maps $e$ also encode a description of
any weighted projective space as the quotient of any other by the
action of a certain 
finite abelian group.

Section \ref{itthis} considers complex oriented cohomology theories
$E^*(-)$, whose coefficient rings $E_*$ are even dimensional and free
of additive torsion. In Theorem \ref{itthompspec}, the complex
orientation is exploited to describe $E^*(\Pc)$ via iterated Thom
isomorphisms and the $E$-theoretic formal group law, whenever $\chi$
is divisive. Illustrative examples are given using integral cohomology
and complex $K$-theory, and the former are shown to recover Kawasaki's
calculations in the divisive (and therefore the $p$-primary) case.

Section \ref{gere} addresses the question that immediately arises: is
it possible to reassemble the $p$-primary parts $\Ppc$, and recover
$\Pc$ up to homeomorphism? Theorem \ref{reass} answers affirmatively,
and offers the surprising addendum that reassembly may be effected in
two contrasting ways, which describe $\Pc$ as an iterated limit and an
iterated colimit of the $\Ppc$s respectively. Both constructions arise
over $\bC P^n$, utilising the maps $e(1/{}_p\chi)$ and $e({}_p\chi)$.
Examples are discussed which show that different reassemblies of the
same $p$-primary parts can produce non-homeomorphic results.

Section \ref{core} introduces cohomological reassembly as a natural
analogue of the previous geometry. Theorem \ref{asscoho} describes
$E^*(\Pc)$ as both a limit and a colimit of the $p$-primary parts
$E^*(\Ppc)$, made explicit in terms of direct sums and iterated tensor
products over $E^*(\bC P^n)$ respectively. Most memorable is the
resulting identification of $E^*(\Pc)$ with an intersection of
$p$-primary subalgebras of $E^*(\bC P^n)$. Examples are provided in
integral cohomology and complex $K$-theory to show that detailed
computation is possible; the former recover Kawasaki's calculations
for arbitrary $\chi$, and the latter those of~\cite{ala:csw}. 

To date, toric topology has rarely dealt with the homology coalgebra
$E_*(X_\varSigma)$ of a toric variety $X_\varSigma$. Section
\ref{hofucl} redresses this situation, at least for weighted
projective space, by studying the universal example
$\varOmega_*^U(\Pc)$ for divisive $\chi$. The relationship between
Bott towers and iterated Thom spaces is recalled, and applied in
Theorem \ref{osuppi} to identify explicit $\osu$-generators for the
coalgebra $\osu(\Pc)$. In particular, a rational fundamental class is
constructed that may be interpreted in terms of toric desingularisation.

Finally, Section \ref{hore} introduces homological reassembly by
dualising the results of Section \ref{core}. Theorem \ref{assho}
describes the coalgebra $E_*(\Pc)$ as both a limit and a colimit of
the $p$-primary parts $E_*(\Ppc)$ in terms of direct sums and iterated
tensor products over $E_*(\bC P^n)$ respectively. In particular,
$E_*(\Pc)$ is identified with an intersection of $p$-primary
subcoalgebras of $E_*(\bC P^n)$, and an illustrative examples is given
in complex $K$-theory.

Throughout our work we write $S^1$ for the circle as a topological
space, and $T<\bC^1$ for its realisation as the group of unimodular
complex numbers with respect to multiplication. For any integer
$k>0$ we write $\bZ/k$ for the integers modulo $k$, and $C_k<T$ for
its realisation as the subgroup generated by a primitive $k$th root
of unity. We interpret the standard simplex $\varDelta^n$ as the
intersection of the positive orthant
$\bR_{\scriptscriptstyle\geqslant}^{n+1}$ with the hyperplane
$x_0+\dots+x_n=1$, and denote its boundary by $\partial\varDelta^n$.
As an abstract simplicial complex, $\partial\varDelta^n$ has
$\binom{n+1}{k+1}$ faces of dimension $k$, for $-1\leq k<n$.

For any generalised cohomology theory, we follow the convention that
all homology and cohomology groups $E_*(X)$ and $E^*(X)$ are
\emph{reduced} for every space $X$. The unreduced counterparts are
given by adjoining a disjoint basepoint, and considering $E_*(X_+)$
and $E^*(X_+)$. The \emph{coefficient ring} $E_*$ is given by
$E_*(S^0) \cong E^{-*}(S^0)$; we identify the homological and
cohomological versions without further comment, and interpret
$E_*(X_+)$ and $E^*(X_+)$ as $E_*$-modules and $E_*$-algebras
respectively. 
We make the important assumption that $E_*$ is even dimensional
and free of additive torsion, as holds for integral cohomology
$H^*(-)$, complex $K$-theory $K^*(-)$, and complex cobordism
$\varOmega^*_U(-)$.

These theories are also \emph{complex oriented} 
\cite[Part II \S 2]{ada:shg}, 
by means of a class $x^E$ in $E^2(\bC P^\infty)$ whose restriction to 
$E^2(\bC P^1)$ is a generator. It follows that there exists a canonical 
isomorphism
\begin{equation}\label{ecp}
E^*(\bC P^\infty_\ssp)\;\cong\;E_*\llbracket x^E\rrbracket
\end{equation}
of $E_*$-algebras, and that complex vector bundles have associated
$E$-theoretic Chern classes. In particular, $x^E$ is the first Chern
class $c_1^E(\zeta)$ of the dual Hopf line bundle $\zeta$ over $\bC
P^\infty$. A minor abuse of notation allows $x^E$ to be confused with
its restriction to $\bC P^n$, and produces an isomorphism
\begin{equation}\label{ecpn} 
E^*(\bC P^n_\ssp)\;\cong\;E_*[x^E]/((x^E)^{n+1})\,.
\end{equation}
It is convenient to denote $x^E$ by $u$ in the universal case
$\varOmega^*_U(\bC P^\infty)$, and to write $x^H$ as $x$ in $H^2(\bC
P^\infty)$.

If the Thom space of $\zeta$ is identified with $\bC P^\infty$, then
$x^E$ may also be interpreted as a Thom class $t^E(\zeta)$, and
extended to a universal Thom class $t^E\in E^0(\MU)$; thus $t^U$ is 
represented by the identity map on $\MU$.

\begin{acknowledgements}
  The authors are grateful to Abdallah Al Amrani, for
  helpful correspondence on the topological history of weighted
  projective space, and for sharing his knowledge of the
  literature. The first and third authors are indebted to Sam Gitler
  and the mathematicians of Cinvestav for their outstanding
  hospitality in Mexico City, and particularly to Ernesto Lupercio for
  proposing the version of Theorem \ref{reass} that appears
  below. The third author thanks Anand Dessai and the University of 
  Fribourg, who provided the opportunity to expound much of this  
  material in May 2008. Credit is also due to Jack Morava for 
  supporting the authors' belief that the topology of weighted 
  projective spaces might well be studied prime by prime.
\end{acknowledgements}

%
%
%
%
%
%
%
%
%

\section{Weighted projective space}\label{wps}

In this section basic notation is established, and the definitions of
weighted projective space, weighted lens space, and their associated
constructions are recalled. Readers are referred to Al Amrani
\cite{ala:csw}, \cite{ala:cbc}, Kawasaki \cite{kaw:ctp}, and the
authors' own work \cite{bafrra:ecr} for further details.

The \emph{standard action} of the $(n+1)$-dimensional torus
$T^{n+1}$ on $\bC^{n+1}$ is by coordinatewise multiplication, and
restricts to the unit sphere $S^{2n+1}$. The orbit space of the
latter is homeomorphic to the standard simplex
$\varDelta^n\subset\bR_{\scriptscriptstyle\geqslant}^{n+1}$, and the
quotient map $r\colon S^{2n+1}\ra\varDelta^n$ is given by
$r(z)=(|z_0|^2,\dots,|z_n|^2)$.

A \emph{weight vector} $\chi$ is a sequence $(\chi_0,\dots,\chi_n)$ of
positive integers, and $T\br{\chi}<T^{n+1}$ denotes the subcircle of
elements $(t^{\chi_0},\dots,t^{\chi_n})$, as $t$ ranges over $T$. It
is convenient to abbreviate the greatest common divisor
$\gcd(\chi_0,\dots,\chi_n)$ and least common multiple
$\lcm(\chi_0,\dots,\chi_n)$ to $g=g(\chi)$ and $l=l(\chi)$
respectively.
\begin{defn}\label{defwproj}
  The \emph{weighted projective space} $\Pc$ is the orbit space of the
  action of $T\br{\chi}$ on $S^{2n+1}$; it admits a \emph{canonical
    action} of the quotient $n$--torus $T^{n+1}/T\langle\chi\rangle$,
  with orbit space $\varDelta^n$.
\end{defn}
The respective quotient maps are
\begin{equation}\label{wpsquot}
S^{2n+1}\stackrel{p(\chi)}{\lllra}
\Pc\stackrel{q(\chi)}{\lllra}\varDelta^n,
\end{equation}
whose composition is $r$. The action of $T\br{\chi}$ is free 
when $\chi=(d,\dots,d)$ for any positive integer $d$, in which case
$\Pc$ reduces to $\CPn$; in general, $\Pc$ has orbifold
singularities. Weighted projective spaces provide an important class
of singular examples in algebraic and symplectic geometry, although
the focus of this article is on their algebraic topology.

It is sometimes convenient to assume that $g(\chi)=1$, because
$T\br{d\chi}$ and $T\br{\chi}$ produce homeomorphic orbit spaces for
all $d$.

\begin{defn}\label{defwlens}
  For any positive integer $k$, the \emph{weighted lens space}
  $L(k;\chi)$ is the orbit space of the action of the weighted cyclic
  subgroup $C_k\br{\chi}<T\br{\chi}$ on $S^{2n+1}$; it admits
  \emph{canonical actions} of the quotient circle
  $T\br{\chi}/C_k\br{\chi}$ with orbit space $\Pc$, and of the
  $(n+1)$--torus $T^{n+1}/C_k(\chi)$ with orbit space $\varDelta^n$.
\end{defn}
If $k$ is prime to $\chi_i$ for $0\leq i\leq n$, then $L(k;\chi)$
is a standard lens space, 
and is smooth; otherwise, $C_k\br{\chi}$ fails to act freely, and 
$L(k;\chi)$ 
may be singular.

Restricting \eqref{wpsquot} to the hyperplane $z_n=0$ yields
orbit maps
\begin{equation}\label{1restrwpsquot}
S^{2n-1}\stackrel{p(\chi')}{\lllra}
\Pd\stackrel{q(\chi')}{\lllra}\varDelta^{n-1},
\end{equation}
where $\chi'$ denotes $(\chi_0,\dots,\chi_{n-1})$ and
$\varDelta^{n-1}$ is the subsimplex $x_n=0$ of $\varDelta^n$.
On the other hand, restricting to the cylinder $|z_n|^2=1/2$ gives
\begin{equation}\label{2restrwpsquot}
S^{2n-1}\times S^1\stackrel{p_{1/2}}{\lllra}
L(\chi_n;\chi')\stackrel{q_{1/2}}{\lllra}\varDelta_{1/2}\,,
\end{equation}
where $S^{2n-1}\subset S^{2n+1}$ is the subsphere
\[
|z_0|^2+\dots+|z_{n-1}|^2\;=\;1/2\,.
\]
Thus $p_{1/2}$ is the orbit map for $T\br{\chi}$, and factors through 
$L(\chi_n;\chi')\times S^1$ under the actions of
$C_{\chi_n}\br{\chi'}$ and $T\br{\chi}/C_{\chi_n}\br{\chi'}$
respectively; the latter is isomorphic to $T$.
Similarly, $q_{1/2}$ is the orbit map for the $n$--torus
$T^n/C_{\chi_n}\br{\chi'}$, whose orbit space
\smash{$\varDelta_{1/2}$} is the $(n-1)$--simplex $x_n=1/2$.  Finally,
restricting to the 
circle $|z_n|=1$ projects every point $(0,\dots,0,z_n)$ in $S^{2n+1}$ 
onto $[0,\dots,0,1]$ in $\Pc$, and thence to the vertex $(0,\dots,0,1)$ 
in $\varDelta^n$.

Now consider the decomposition of $\varDelta^n$ into the union of
subspaces $N(1/2)$ and $C(1/2)$, specified by $x_n\leq 1/2$ and
$x_n\geq 1/2$ respectively; they are homeomorphic to the product
$\varDelta^{n-1}\times [0,1/2]$ and the cone
\smash{$C\varDelta_{1/2}$}. So $\Pc$ may be expressed as the
pushout of
\begin{equation}\label{pchicolim}
N\Pd\stackrel{i}{\llla}L(\chi_n;\chi')
\stackrel{j}{\llra}CL(\chi_n;\chi'),
\end{equation}
where $N\Pd$ denotes the neighbourhood $q^{-1}(N(1/2))$ of $\Pd$ in
$\Pc$, and $CL(\chi_n;\chi')$ denotes the cone $q^{-1}(C(1/2))$,
with basepoint $[0,\dots,0,1]$. Equivalently, \eqref{pchicolim}
arises by decomposing $S^{2n+1}$ as the pushout of
\begin{equation*}\label{sphcolim}
S^{2n-1}\times D^2\stackrel{i}\llla S^{2n-1}\times S^1
\stackrel{j}{\llra}D^{2n}\times S^1,
\end{equation*}
and forming orbit spaces under the action of $T\br{\chi}$.
Reparametrising $D^{2n}$ shows that $[0,\dots,0,1]$ admits a
neighbourhood of the form $\bC^n\big/C_{\chi_n}$; repeating at
each point $[0,\dots,0,1,0,\dots,0]$ confirms that $\Pc$ is a complex
orbifold.

Diagram \eqref{pchicolim} is cofibrant, and therefore expresses
$\Pc$ as the homotopy colimit of the diagram
\begin{equation}\label{pchiishocolim}
\Pd\stackrel{f}{\llla}L(\chi_n;\chi')\llra *\:,
\end{equation}
where $f$ denotes the orbit map for the circle
$T\br{\chi'}/C_{\chi_n}\br{\chi'}$, and $*$ is the point
$[0,\dots,0,1]$. This reinterprets Kawasaki's cofibre sequence
\cite[page 245]{kaw:ctp}
\begin{equation}\label{cofiblww}
L(\chi_n;\chi')\stackrel{f}{\llra}\Pd\stackrel{g}{\llra}\Pc.
\end{equation}
\begin{rem}\label{betterhocolim}
  The category underlying diagram \eqref{pchiishocolim} may also be
  construed as $\cat{cat}(\partial\varDelta^1)$, whose objects are the
  faces $\varnothing$, $0$ and $1$ of $\partial\varDelta^1$ and
  morphisms their inclusions. Iteration on $\Pd$ leads to a
  description of $\Pc$ as a homotopy colimit over
  $\cat{cat}(\partial\varDelta^n)$, in which the relevant diagram
  assigns an orbit space $T^n/T^k(\sigma)$ to each
  $(k-1)$-dimensional face $\sigma$ of $\varDelta^n$; this is
  precisely the homotopy colimit of \cite[Proposition
  5.3]{wezizi:hcc}.
\end{rem}

Following Kawasaki, Al Amrani \cite[I.1(b)]{ala:csw} defines maps
\begin{equation}\label{defphipsi}
\CPn\stackrel{\phi}{\llra}\Pc\stackrel{\psi}{\llra}\CPn
\end{equation}
by $\phi[z_0,\dots,z_n]=[z_0^{\chi_0},\dots,z_n^{\chi_n}]$ and
$\psi[z_0,\dots,z_n]=
[z_0^{l(\chi)/\chi_0},\dots,z_n^{l(\chi)/\chi_n}]$. In both cases, the
formulae for the homogeneous coordinates of the target values are
understood to be normalised. It is sometimes important to make
the weights explicit, by writing $\phi(\chi)$ and $\psi(\chi)$
respectively.

Usually, $\psi$ is interpreted as a complex line bundle over $\Pc$,
but may equally well be specified by its first Chern class $c_1(\psi)$
in $H^2(\Pc;\bZ)$. The composition $\psi\scirc\phi$ has degree 
$l=l(\chi)$ on $H^2(\CPn;\bZ)$, so $\phi^*(\psi)$ is the $l$th tensor 
power $\zeta^l$ of the dual Hopf bundle, and $c_1(\phi^*(\psi))=lx$ in 
$H^2(\CPn)$, following \eqref{ecp}. 

The r\^ole of $\psi$ is clarified by identifying the total space
$S(\psi)$ of its associated circle bundle. 
\begin{prop}\label{llchi}
The space $S(\psi)$ is a $(2n+1)$-dimensional weighted lens space 
$L(l;\chi)$.
\end{prop}
\begin{proof}
By definition, $S(\psi)$ is the pullback of the diagram
\[
\Pc\stackrel{\psi}{\llra}\bC P^n\llla S^{2n+1},
\]
and is a subspace $X\subset S^{2n+1}\times\Pc$. It contains all pairs
$(y,[z])$ that satisfy the equation
\begin{equation}\label{pairs}
t\pfm(y_0,\dots,y_n)\;=\;
(z_0^{l(\chi)/\chi_0},\dots,z_n^{l(\chi)/\chi_n})
\end{equation}
in $S^{2n+1}$, for some $t\in T$. So there exists a map $h\colon
L(l(\chi);\chi)\ra X$, defined by
\[
h[w_0,\dots,w_n]\;=\;
\bigl((w_0^{l(\chi)/\chi_0},\dots,w_n^{l(\chi)/\chi_n}),
[w_0,\dots,w_n]\bigr)\,;
\]
moreover, an inverse to $h$ is given by mapping $(y,[z])$ to the
equivalence class of those $(n+1)$-tuples $(z_0,\dots,z_n)$ for which
$t=1$ in \eqref{pairs}. It follows that $h$ is the required
homeomorphism.
\end{proof}
\begin{cor}\label{pcisfreequot}
The circle $T\br{\chi}/C_l\br{\chi}$ acts freely on $L(l;\chi)$, and
has orbit space $\Pc$.
\end{cor}

Of course the associated sphere bundle $S(\phi^*(\psi))$ is
homeomorphic to the standard lens space $L(l;1,\dots,1)$, and is
therefore a smooth manifold.

The following natural numbers are associated to $\chi$, and were
essentially introduced by Kawasaki; alternative descriptions are
recovered in Theorem \ref{kawa}.
\begin{defn}\label{deflkqk}
  For any subset $J\subseteq[n]$, the integer $\chi_J$ is the product
  $\prod_{j\in J}\chi_j$, and $h_J=h_J(\chi)$ is the quotient
  $\chi_J/\gcd(\chi_j:j\in J)$; for any $1\leq j\leq n$, the integer
  $l_j=l_j(\chi)$ is $\lcm(h_J:|J|=j)$, and $m_j=m_j(\chi)$ is
  $l(\chi)^j/l_j$.
\end{defn}
Thus $l_1=l$ and $m_1=1$, whereas
$l_n=\chi_0\cdots\chi_n/g$ and $m_n=g\3 l^n/\chi_0\cdots\chi_n$.

Kawasaki applies the cofibre sequence \eqref{cofiblww} to identify the
integral cohomology ring of $\Pc$ by means of an isomorphism
\begin{equation}\label{cohringwps}
H^*(\Pc_\ssp;\bZ)\;\cong\;\bZ[v_1,\dots,v_n]/I(\chi)\,,
\end{equation}
where $v_j$ has dimension $2j$ and $I(\chi)$ is the ideal generated
by the elements \smash{$v_1^j-m_jv_j$} for $1\leq j\leq n$. Moreover, $v_1$
equals $c_1(\psi)$, so $\phi^*(v_1)=lx$ holds in $H^2(\CPn;\bZ)$. The
same calculations identify the integral cohomology ring of
$L(\chi_n;\chi')$ in terms of additive isomorphisms
\begin{equation}\label{cohl}
H^j(L(\chi_n;\chi');\bZ)\;\cong\;
\begin{cases}
\bZ&\text{if $j=2n-1$}\\
\bZ/s_k&\text{if $j=2k$\; for\; $1\leq k\leq n-1$}\\
0&\text{otherwise,}
\end{cases}
\end{equation}
where $s_k=s_k(\chi)$ is given by $l_k/l'_k$, and $l'_k\letbe l_k(\chi')$.
\begin{rem}\label{rattype}
  These results reveal that the maps $\phi$ and $\psi$ of
  \eqref{defphipsi} are mutually inverse rational homotopy
  equivalences, as their definitions suggest. The rationalisation
  $\Pc_\bQ$ is therefore the $2n$--skeleton of an Eilenberg--Mac Lane
  space $K(\bQ,2)$.
\end{rem}

Many of Kawasaki's calculations are recovered in Theorem \ref{kawa}.

%
%
%
%
%
%
%
%
%

\ifnige
\section{Iterated Thom spaces and $p$-primary parts}
\else
\section{Iterated Thom spaces and \texorpdfstring{$p$}{p}-primary parts}
\fi
\label{itthspprpa}

In order to follow homotopy theoretical convention and study $\Pc$ 
prime by prime, it is convenient to introduce certain restrictions on
the weights. 

If $g(\chi)=1$, then there exists an isomorphism
\[
\bP(d\chi_0,\dots,d\chi_{j-1},\chi_j,d\chi_{j+1},\dots,d\chi_n)
\;\cong\;\Pc
\]
of 
algebraic varieties for any natural number $d$ such that
$\gcd(d,\chi_j)=1$, and every $0\leq j\leq n$ \cite{dol:wpv}, 
\cite{ian:wwc}. So no generality is lost by insisting that $\chi$ is
\emph{normalised}, in the sense that
\begin{equation}\label{norm}
\gcd(\chi_0,\dots,\widehat{\chi_j},\dots,\chi_n)\;=\;1
\end{equation}
for every $0\leq j\leq n$.
\begin{defns}\label{fidivtodiv}
  The weight vector $\chi$ is
\begin{enumerate}
\item[\bf 1.] 
\emph{weakly divisive} if $\chi_j$ divides $\chi_n$ for every $0\leq j<n$
\item[\bf 2.] 
\emph{divisive} if $\chi_{j-1}$ divides $\chi_j$ for every $1\leq j\leq n$
\item[\bf 3.] 
\emph{$p$-primary} if every $\chi_j$ is a power of a fixed prime $p$.
\end{enumerate}
\noindent
If $\chi$ is divisive, then $q_j=q_j(\chi)$ denotes the integer 
$\chi_j/\chi_{j-1}$, for $1\leq j\leq n$.
\end{defns}
So divisive implies weakly divisive, but not conversely. In fact $\chi$ 
is divisive precisely when the reverse sequence $\chi_n$, \dots, $\chi_0$ 
is \emph{well ordered}, in the sense of Nishimura--Yosimura 
\cite{niyo:qkt}; then 
\begin{equation*}
  * = \bP(\chi_{n}),\;\;
  \bP(\chi_{n-1},\chi_{n})\setminus \bP(\chi_{n}),\;\;
  \dots,\;\;
  \bP(\chi_{0},\dots,\chi_{n})\setminus \bP(\chi_{1},\dots,\chi_{n})
\end{equation*}
is a cell decomposition of $\Pc$ with one cell in every even dimension
(compare \cite[Remark 3.2]{bafrnora:cwp} and 
\cite[Proposition 2.3]{niyo:qkt}). This decomposition describes the 
canonical cells that arise from Corollary \ref{divisprop} below. 
\begin{lem}\label{pspeclem}
Every normalised $p$-primary weight vector $\pi$ may be ordered so as 
to take the form
\begin{equation}\label{pspec}
\big(1,1,p^{k_2},\dots,p^{k_n}\,\big)
\end{equation}
for some sequence $k_2\leq\dots\leq k_n$ of exponents.
\end{lem}
\begin{proof}
Since $g(\pi)=1$, it follows that $\pi_j=1$ for some $j$. Applying 
\eqref{norm}, and reordering if necessary, yields the required form.
\end{proof}
Of course \eqref{pspec} is divisive, and identifies $k_0$ and $k_1$ as 
$0$.

\begin{defn}\label{defitthomcplx}
  If $(X_i)$ is a sequence of topological spaces for $0\leq i\leq n$,
  and $(\theta_i)$ a sequence of vector bundles over $X_i$
  for $0\leq i<n$, then $X_n$ is an \emph{$n$-fold iterated Thom
    space} over $X_0$ whenever $X_i$ is homeomorphic to the Thom
  space $\Th(\theta_{i-1})$ for every $1\leq i\leq n$.
\end{defn}
\begin{exa}\label{cpnisitc}
  If $\zeta_i$ denotes the dual Hopf line bundle over $\bC P^i$ for
  every $i\geq 0$, then the standard homeomorphisms $\bC
  P^i\cong\Th(\zeta_{i-1})$ display $\CPn$ as an $n$-fold iterated Thom
  space over the one-point space $\bC P^0$; or, alternatively, as an
  $(n-1)$-fold iterated Thom space over the $2$-sphere $\bC P^1$.
\end{exa}
In Example \ref{cpnisitc}, $\zeta_i$ may equally well be replaced by
the Hopf bundle itself.

\begin{thm}\label{wdivisprop}
  If $\chi$ is weakly divisive, then $\Pc$ is homeomorphic to 
the Thom space of a complex line bundle over $\bP(\chi')$.
\end{thm}
\begin{proof}
  Al Amrani's proof of \cite[I.1(c)]{ala:csw} applies to the
  line bundle $\psi'$ over $\bP (\chi')$, and shows that the unit
  disc bundle $D(\psi')$ is homeomorphic to the neighbourhood $N\Pd$
  of $\Pd\subset\Pc$, as defined in \eqref{pchicolim}. The unit sphere
  bundle $S(\psi')$ is therefore homeomorphic to the weighted lens
  space $L(\chi_n;\chi')$, and the cofibre sequence \eqref{cofiblww}
  identifies $\Pc$ with the Thom space $\Th(\psi')$.
\end{proof}
\begin{cor}\label{divisprop}
  If $\chi$ is divisive, then $\Pc$ is homeomorphic to an $n$-fold
  iterated Thom space of complex line bundles over the one-point
  space $*$\,.
\end{cor} 

By analogy with Remark \ref{betterhocolim}, an $n$-fold iterated Thom
space may also be expressed as an iterated pushout, and therefore as
a homotopy colimit over the category $\cat{cat}(\partial\varDelta^n)$.

In order to apply Corollary \ref{divisprop} further, the maps $\phi$
and $\psi$ of \eqref{defphipsi} must be generalised. A suitable
context is provided by interpreting weight vectors as elements of the
multiplicative monoid $\bN^{n+1}$, with identity element
$1=(1,\dots,1)$.  Given any two such $\chi$ and $\omega$, there exists
a smallest positive integer $s=s(\chi,\omega)$ such that $\omega$
divides $s\chi$. The resulting quotient has coordinates
$s\chi_j/\omega_j$ for $0\leq j\leq n$, and is conveniently denoted by
$\chi/\omega$; so equations such as
\begin{equation}\label{wteqns}
\begin{split}
s\;=\;(s,&\dots,s),\quad\omega(\chi/\omega)\;=\;s\chi,
\quad(\omega\chi)/\omega\;=\;\chi,\\
&\chi\;=\;\chi/1,\sands 1/\chi\;=\;l(\chi)/\chi
\end{split}
\end{equation}
hold amongst weight vectors.

Every $\chi$ may then be expressed as a product of indecomposables. For 
any $0\leq j\leq n$ and any prime $p$, write $\chi_j$ as 
$p^{a(j)}\alpha_j$, where $p^{a(j)}$ denotes the $p$-content of $\chi_j$ 
and $\gcd(p,\alpha_j)=1$.
\begin{defn}\label{pprim}
The \emph{$p$-content} of~$\chi$ is the $p$-primary weight vector
\[
{}_{p}\chi\;\letbe\;(p^{a(0)},\dots,p^{a(n)})\,,
\]
which satisfies $\chi={}_p\chi\alpha$ in $\bN^{n+1}$; the
\emph{primary decomposition} of $\chi$ is the factorisation
$\chi={}_{p_1}\chi\cdots{}_{p_m}\chi$, as $p_i$ ranges over the
prime factors of the $\chi_j$.
\end{defn}
If $\chi$ is normalised then so is ${}_p\chi$, but the non-decreasing
property is \emph{not} hereditary in this sense. It follows from
Definition \eqref{pprim} that $\alpha=\prod_{p_i\neq p}{}_{p_i}\chi$,
and that \smash{$l(\chi)=p^{m(a)}l(\alpha)$} where
$m(a)\letbe\max_ia(i)$.

\begin{rem}\label{htpytype}
   Recent results of \cite{bafrnora:cwp} show that, amongst weighted
   projective spaces, 
   the homotopy type of $\Pc$ is determined by the \emph{unordered}
   coordinates of its non-trivial $p$-contents ${}_p\chi$, for normalised
   $\chi$. The ${}_p\chi$ may therefore assumed to be non-decreasing, and
   remultiplied to give a weight vector $\chi^*$ for which there exists
   a homotopy equivalence $\Pc\simeq\bP(\chi^*)$. Moreover $\chi^*$
   is divisive by construction, so Corollary \ref{divisprop} implies
   that every $\Pc$ is homotopy equivalent to an iterated Thom space.
   If the weights are pairwise coprime then $\chi^*$ takes the form
   $(1,\dots,1,c)$ and $\bP(\chi^*)$ reduces to $\Th(\zeta^c_{n-1})$ over 
   $\bC P^{n-1}$, where $c=\prod_i\chi_i$.
\end{rem}
\begin{defns}\label{defec}\hfill
\begin{enumerate}
\item[\bf 1.]
The map $e(\chi/\omega)\colon\Po\ra\Pc$ is given by
\[
e(\chi/\omega)([z_0,\dots,z_n])\;=\;
[z_0^{s\chi_0/\omega_0},\dots,z_n^{s\chi_n/\omega_n}]\,,
\]
where $s=s(\chi,\omega)$, and coordinates are normalised as necessary:
\item[\bf 2.] 
the group $C_{\chi/\omega}$ is the product 
\[
C_{s\chi_0/\omega_0}\times\dots\times C_{s\chi_n/\omega_n}
\]
of cyclic groups, considered as a subgroup of $T^{n+1}$.
\end{enumerate}
\end{defns}
Following \eqref{wteqns}, the cases $e(\phi\chi/\chi)$ and
$C_{\phi\chi/\chi}$ reduce to $e(\phi)$ and $C_\phi$ respectively. By
definition, $e(r)=e(r,\dots,r)$ raises homogeneous coordinates in
$\Pc$ to the $r$th power, and is therefore known as the \emph{$r$th
  power map} on $\Pc$.
\begin{prop}\label{eisfg}
The map $e(\chi/\omega)$ is the orbit map of the natural action of 
$C_{\chi/\omega}$ on $\Po$.
\end{prop} 
\begin{proof}
  Note first that
  $e(\chi/\omega)([y_0,\dots,y_n])=e(\chi/\omega)([z_0,\dots,z_n])$
  holds in $\Pc$ precisely when
\[
[y_0,\dots,y_n]\;=\;
\{[\lambda_0z_0,\dots,\lambda_nz_n]:\lambda_0^{h\omega_0/\chi_0}=\dots
=\lambda_n^{h\omega_n/\chi_n}=1\}
\] 
  in $\Po$. Since $e(\chi/\omega)$ is clearly surjective, the result
  follows.
\end{proof}
\begin{cor}
Any weighted projective space arises as the orbit space of any other of 
the same dimension, under the action of a finite abelian group. 
\end{cor}
\begin{rem}\label{erems}
  In the language of Definition \ref{defec} and Proposition
  \ref{eisfg}, Al Amrani's maps $\phi(\chi)$ and $\psi(\chi)$ are
  given by $e(\chi)$ and $e(1/\chi)$ respectively.  Kawasaki
  \cite[page 243]{kaw:ctp} notes that $\phi(\chi)$ is the orbit map of
  the action of $C_\chi$ on $\CPn$.
\end{rem}

\begin{prop}\label{ecomp}
For any weight vectors $\omega$ and $\chi$, the composition
\[
\Po\stackrel{e(\sigma/\omega)}{\llllra}\bP(\sigma)
\stackrel{e(\chi/\sigma)}{\llllra}\Pc
\]
factorises as $e(s')\scirc e(\chi/\omega)=e(\chi/\omega)\scirc e(s')$,
where $s'$ denotes the natural number
$s(\omega,\sigma)s(\sigma,\chi)/s(\omega,\chi)$.
\end{prop}
\begin{proof}
  It suffices to note that the given composition acts on homogeneous
  coordinates in $\Pc$ by $z_i\mapsto
  z_i^{s(\omega,\sigma)s(\sigma,\chi)\chi_i/\omega_i}$, for $0\leq i\leq n$.
\end{proof}
Proposition \ref{ecomp} implies the factorisations  
\begin{equation*}
e(\chi/\omega)\scirc e(\omega/\chi)=e(s'')\sands
e(\omega/\chi)\scirc e(\chi/\omega)=e(s'')\,, 
\end{equation*}
where $s''\letbe s(\omega,\chi)s(\chi,\omega)$. Similarly,
$e(\chi)\scirc e(1/\omega)=e(l(\omega)/s(\omega,\chi))$
for any weight vectors $\chi$ and $\omega$.

\begin{defn}\label{defees}
  For any weight vector $\chi$ and any prime $p$, the
  \emph{$p$-primary part} of $\Pc$ is the weighted projective space
  $\Ppc$; the canonical maps 
\[
e({}_p\chi/\chi)\colon\Pc\lra\Ppc\sands e(\chi/{}_p\chi)\colon\Ppc\lra\Pc
\]
are \emph{$p$-extraction} and \emph{$p$-insertion} respectively.
\end{defn}
In the notation of Definition \ref{pprim}, extraction and insertion 
are given by
\begin{equation}\label{epcoord}
\begin{split} e({}_p\chi/\chi)[z_0,\dots,z_n]\;=\;
[z_0^{l(\alpha)/\alpha_0},&\dots,z_n^{l(\alpha)/\alpha_n}]
\\\text{and}\quad &e(\chi/{}_p\chi)[z_0,\dots,z_n]\;=\;
[z_0^{\alpha_0},\dots,z_0^{\alpha_n}]\,,
\end{split}
\end{equation}
in terms of homogeneous coordinates. By Proposition \ref{ecomp}, the
compositions $e(\chi/{}_p\chi)\scirc e({}_p\chi/\chi)$ and
$e({}_p\chi/\chi)\scirc e(\chi/{}_p\chi)$ reduce to the appropriate
power maps $e(l(\alpha))$. 
\begin{rem}
It follows immediately from Lemma \ref{pspeclem} and Corollary 
\ref{divisprop} that every $p$-primary part $\Ppc$ is an iterated Thom 
space over $*$.
\end{rem}

\begin{exa}\label{345a}
  The $2$-, $3$-, and $5$-primary parts of $\bP(3,4,5)$ are $\bP(1,4,1)$,
  $\bP(3,1,1)$, and $\bP(1,1,5)$ respectively. They form the codomains of 
  the $2$-, $3$-, and $5$-extraction maps, whose values on $[z_0,z_1,z_2]$ 
  are given by  
\[
[z_0^5,z_1^{15},z_2^3],\quad [z_0^{20},z_1^5,z_2^4],\sands
[z_0^4,z_1^3,z_2^{12}]
\]
respectively. The $2$-, $3$-, and $5$-primary parts also form the
domains of the $2$-, $3$-, and $5$-insertion maps, whose values on
$[z_0,z_1,z_2]$ are given by
\[
[z_0^3,z_1,z_2^5],\quad [z_0,z_1^4,z_2^5],\sands
[z_0^3,z_1^4,z_2]
\]
respectively.
\end{exa}

%
%
%
%
%
%
%
%
%

\section{Geometric reassembly}\label{gere}

The problem of reassembling $\Pc$ from its $p$-primary parts must now
be addressed. The solution is best understood in terms of weight vectors 
$\sigma$ and $\sigma'$, and commutative squares of the form
\begin{equation}\label{chirhos}
\begin{CD}
\bP(\sigma\sigma')@>{e(1/\sigma')}>>\bP(\sigma)\\
@V{e(1/\sigma)}VV@VV{e(1/\sigma)}V\\
\bP(\sigma')@>>{e(1/\sigma')}> 
\bC P^n
\end{CD}
\quad\spandsp\quad
\begin{CD}
\bC P^n@>{e(\sigma)}>>\bP(\sigma)\\
@V{e(\sigma')}VV@VV{e(\sigma')}V\\
\bP(\sigma')@>>{e(\sigma)}> 
\bP(\sigma\sigma')
\end{CD}
\quad;
\end{equation}
these may be incorporated into a single generalised square, although
restrictions must be imposed upon $\sigma$ and $\sigma'$.

For any weight vector $\chi$, it is convenient to write
$Q(\chi)\subset\bZ$ for the set of primes that occur non-trivially in
its $p$-primary decomposition.
\begin{defn}\label{wcc}
  Two weight vectors $\sigma$ and $\sigma'$ are \emph{coprime} if
  $Q(\sigma)\cap Q(\sigma')=\varnothing$; in other words, if
  $\gcd(\sigma_j,\sigma'_k)=1$ for every $0\leq j,k\leq n$.
\end{defn}
In the context of Definition \ref{defec}, this condition is equivalent
to the coprimality of $|C_\sigma|$ and $|C_{\sigma'}|$. It also
implies the weaker condition, that
\begin{equation}\label{chineseremainder}
h(g_0,\dots,g_n)\;=\;
\big((g_0^{\sigma_0},\dots,g_n^{\sigma_n}),
(g_0^{\sigma'_0},\dots,g_n^{\sigma'_n})\big)
\end{equation}
defines an isomorphism $h\colon C_{\sigma\sigma'}\to C_{\sigma'}\times
C_{\sigma}$; or equivalently, that $C_{\sigma\sigma'}$ is generated by the
subgroups $C_{\sigma}$ and $C_{\sigma'}$. 
\begin{rem}\label{qsqsinv}
Since $Q(\sigma)=Q(1/\sigma)$, it follows that $\sigma$ and $\sigma'$ 
are coprime if and only if $1/\sigma$ and $1/\sigma'$ are coprime.
\end{rem} 
\begin{prop}\label{pbs}
  If $\sigma$ and $\sigma'$ are coprime, then the left-hand diagram 
  {\rm \ref{chirhos}} is a pullback square.
\end{prop}
\begin{proof}
  The pullback of the diagram $\bP(\sigma')\ra\bC P^n\la\bP(\sigma)$ 
  is the subspace
\[
P\;\letbe\;\{(z',z):e(1/\sigma')(z')=e(1/\sigma)(z)\}\;\subset\;
\bP(\sigma')\times\bP(\sigma)\,,
\]
and the canonical map $f\colon\bP(\sigma\sigma')\to P$ acts by 
$f(z)=(e(1/\sigma)(z),e(1/\sigma')(z))$. Since $P$ is Hausdorff, it 
$f$ is surjective, by equivariance.
remains to show that $f$ is bijective.

The maps $e(1/\sigma)$ and $e(1/\sigma')$ are quotients by $C_{1/\sigma}$ 
and $C_{1/\sigma'}$ respectively, by Proposition \ref{eisfg}. So
$C_{1/\sigma'}\times C_{1/\sigma}$ acts on $P$, and $f$ is equivariant 
with respect to the isomorphism $h$ of \eqref{chineseremainder}. The
projection $p\colon P\to\bP(\sigma)$ is the quotient by
$C_{1/\sigma}$, and is equivariant with respect to the projection
$C_{1/\sigma'}\times C_{1/\sigma}\to C_{1/\sigma'}$; the corresponding
statement for $p'\colon P\to\bP(\sigma')$ holds by analogy. Both
$e(1/\sigma)$ and $p$ are surjective, and the former factorises as
$p\scirc f\colon\bP(\sigma\sigma')\to\bP(\sigma)$; therefore 
$f$ is surjective, by equivariance.

To confirm injectivity, choose $x\in\bP(\sigma\sigma')$ and let 
$G<C_{1/\sigma}$ be the isotropy group of $y=f(x)$. Then 
$c\letbe |f^{-1}(y)|=|Gx|$ divides $|G|$, so $c$ divides $|C_{1/\sigma}|$. 
Replacing $\sigma$ by $\sigma'$ and applying the  
corresponding reasoning shows that $c$ also divides $|C_{1/\sigma'}|$. But 
$\sigma$ and $\sigma'$ are coprime, so $c=1$ as sought.
\end{proof}
\begin{rem}
  Since all the maps in the square are algebraic, Proposition
  \ref{pbs} also holds in the category of complex algebraic varieties.
\end{rem}

The following example shows that the square is not generally 
a pullback. 
\begin{exa}\label{nonpb}
If $\sigma=(1,2,2)$ and $\sigma'=(2,1,2)$, then 
$\bP(\sigma)\cong\bP(\sigma')\cong\bC P^2$ and $e(1/\sigma)$ and 
$e(1/\sigma')$ are homeomorphisms. So the pullback is also homeomorphic 
to $\bC P^2$, and cannot be $\bP(\sigma\sigma')\cong\bP(1,1,2)$ 
because the latter contains a singular point. Proposition \ref{eisfg}
confirms that the canonical map $f\colon\bP(\sigma\sigma')\to\bC P^2$ is 
the quotient by $C_2\times C_2$.

An example with normalised weights is given by increasing the 
dimensions, with additional weights $1$; thus $\sigma=(1,1,2,2)$ and 
$\sigma'=(1,2,1,2)$.
\end{exa}

Remarkably, weighted projective spaces may also be expressed as pushouts.
\begin{prop}\label{pos}
  If $\sigma$ and $\sigma'$ are coprime, then the right-hand diagram 
  {\rm \ref{chirhos}} is a pushout square.
\end{prop}
\begin{proof}
The pushout of the diagram $\bP(\sigma')\la\bC P^n\ra\bP(\sigma)$ 
is the space 
\[
R\;=\;\big(\bP(\sigma')\pfm\sqcup\,\bP(\sigma)\big)\big/\!\sim\,,
\]
where the equivalence relation is generated by 
$e(\sigma')(z)\sim e(\sigma)(z)$ for any $z$ in $\bC P^n$; the canonical 
map $g\colon R\to\bP(\sigma\sigma')$ is defined by 
\[
g([e(\sigma')(z)])\;=\;g([e(\sigma)(z)])\;=\;e(\sigma\sigma')(z)\,.
\]
Since $R$ is compact, it suffices to show that $g$ is bijective.

Both $e(\sigma)$ and $e(\sigma')$ are surjective, and the natural map
$q\colon\bC P^n\to R$ is given by
$q(z)=[e(\sigma)(z)]=[e(\sigma')(z)]$; so $q$ is surjective. It is
also invariant with respect to the action of $C_{\sigma\sigma'}$,
since the latter is generated by the subgroups $C_{\sigma'}$ and
$C_{\sigma}$ (via \eqref{chineseremainder}), which act trivially on
$\bP(\sigma)$ and $\bP(\sigma')$ respectively. Hence $q$ induces a
surjective map $\overline{q}\colon\bP(\sigma'\sigma)\to R$.
 
The map $e(\sigma'\sigma)\colon\bC P^n\to\bP(\sigma'\sigma)$
factorises as $g\scirc q$, so $g\scirc\overline{q}=1$ on
$\bP(\sigma'\sigma)$, and $\overline{q}$ is injective. Thus
$\overline{q}$ is bijective, and has inverse $g$, as required.
\end{proof}

In order to state the main reassembly theorem, it is convenient to write 
$[m]$ for the simplicial complex consisting of $m$ disjoint vertices.
\begin{thm}\label{reass}
The weighted projective space $\Pc$ is homeomorphic to the limit of
the $\cat{cat}^{op}[m]$-diagram
\[
\raisebox{4pt}{$\Ddots$}\hspace{-60pt}
\begin{CD}
@.\bP({}_{p_i}\chi)@.\\
@.@VV{e(1/{}_{p_i}\chi)}V@.\\
\bP({}_{p_1}\chi)@>>{e(1/{}_{p_1}\chi)}>\CPn@<<{e(1/{}_{p_m}\chi)}<
\bP({}_{p_m}\chi)\,,
\end{CD}
\hspace{-66pt}\raisebox{4pt}{$\ddots$}\\
\]
and the associated universal maps $\Pc\ra\bP({}_{p_i}\chi)$ may be 
identified with $e({}_{p_i}\chi\?/\chi)$ for every $0\leq i\leq m$.  
Similarly, $\Pc$ is also homeomorphic to the colimit of the 
$\cat{cat}[m]$-diagram
\[
\raisebox{4pt}{$\Ddots$}\hspace{-56pt}
\begin{CD}
@.\bP({}_{p_i}\chi)@.\\
@.@AA{e({}_{p_i}\chi)}A@.\\
\bP({}_{p_1}\chi)@<<{e({}_{p_1}\chi)}<\CPn@>>{e({}_{p_m}\chi)}>
\bP({}_{p_m}\chi)\,,
\end{CD}
\hspace{-62pt}\raisebox{4pt}{$\ddots$}\\
\]
and the associated universal maps $\bP({}_{p_i}\chi)\ra\Pc$ may
be identified with $e(\chi/{}_{p_i}\?\chi)$ for every $0\leq i\leq m$.
\end{thm}
\begin{proof}
Proceed by induction on $m$, noting that the results are trivial for
$m=1$.

Suppose that $Q(\chi)=\{p_1,\dots,p_k,p\}$ and
$\chi_i=p^{a(i)}\alpha_i$, as in Definition \ref{pprim}; so
$Q(\alpha)=\{p_1,\dots,p_k\}$. By the inductive hypotheses,
$\bP(\alpha)$ is homeomorphic to the pullback of the
$\bP({}_{p_i}\chi)$ along the maps $e(1/{}_{p_1}\chi)$, and to the
pushout of the $\bP({}_{p_i}\chi)$ along the maps $e({}_{p_1}\chi)$;
also, the universal maps $\bP(\alpha)\ra\bP({}_{p_i}\chi)$ and
$\bP({}_{p_i}\chi)\ra\bP(\alpha)$ are given by $e({}_{p_i}\chi/\chi)$
and $e(\chi/{}_{p_i}\chi)$ respectively, for $1\leq i\leq k$. It
therefore remains to prove that $\Pc$ is homeomorphic to the pullback
of
\begin{equation}\label{inducpb}
\Ppc\stackrel{e(1/{}_p\chi)}{\llllra}\CPn
\stackrel{e(1/\alpha)\!}{\llllla}\bP(\alpha)
\end{equation}
and the pushout of 
\begin{equation}\label{inducpo}
\Ppc\stackrel{e({}_p\chi)}{\llllla}\CPn
\stackrel{e(\alpha)\!}{\llllra}\bP(\alpha)
\end{equation}
respectively, and to confirm the identity of the associated maps
$\Pc\ra\Ppc$, $\Pc\ra\bP(\alpha)$, $\Ppc\ra\Pc$, and
$\bP(\alpha)\ra\Pc$. These follow directly from Proposition \ref{pbs}
and Proposition \ref{pos} respectively, because ${}_p\chi$ and
$\alpha$ are coprime. The induction is then complete.
\end{proof}
\begin{exa}\label{345b}
  Theorem \ref{reass} applies to Example \ref{345a}, and expresses 
  $\bP(3,4,5)$ as the limit and colimit of the $\cat{cat}[3]$-diagrams
\[
\begin{CD}
@.\bP_3@.\\
@.@VVe_3V@.\\
\bP_2@>>e_2>\bC P^2@<<e_5<\bP_5
\end{CD}\spandsp
\begin{CD}
@.\bP_3@.\\
@.@AAe_3A@.\\
\bP_2@<<e_2<\bC P^2@>>e_5>\bP_5
\end{CD}
\]
respectively, where $\bP_2\letbe\bP(1,4,1)$, $\bP_3\letbe\bP(3,1,1)$, 
and $\bP_5\letbe\bP(1,1,5)$.
\end{exa}
Observe that $\bP(1,5,12)$, $\bP(1,4,15)$, $\bP(1,3,20)$, and
$\bP(1,1,60)$ may also be obtained by recombining $\bP(1,4,1)$,
$\bP(3,1,1)$, and $\bP(1,1,5)$ with permuted coordinates. No two of
the four are homeomorphic, as their singularity structure shows; but
the results of \cite{bafrnora:cwp} (as described in Remark
\ref{htpytype} above) prove that all four are homotopy equivalent to
$\bP(3,4,5)$. The fact that their cohomology rings are isomorphic is
noted in \cite{bafrra:ecr}, and reproven in Theorem \ref{kawa} below;
it is also, of course, implicit in \cite{kaw:ctp}.

%
%
%
%
%
%
%
%
%

\section{Iterated Thom isomorphisms}\label{itthis}

From this point onwards, $E^*(-)$ denotes a complex oriented
cohomology theory, with orientation class $x^E$. As described in
Section \ref{intro}, the crucial examples are: $H^*(-)$, with the Thom
orientation; $K^*(-)$, with the Conner-Floyd orientation; and
$\ous(-)$, with the universal orientation. In particular, $H^*(X)$
denotes the reduced integral cohomology ring of any space $X$.

The existence of a Thom class $t^E$ leads to the \emph{Thom isomorphism}, 
which features in the following standard result. 
\begin{prop}\label{thomstr}
  For any $k$-dimensional complex vector bundle $\theta$ over $X$, the
  $E_*$-algebra $E^*(\Th(\theta))$ is a free module over $E^*(X_\ssp)$ on
  the single generator $t^E(\theta)$; its multiplicative structure is
  determined by the relation
\begin{equation}\label{multreln1}
(t^E(\theta))^2\;=\;c^E_k(\theta)\cdot t^E(\theta)
\end{equation}
in $E^{4k}(\Th(\theta))$.
\end{prop}
More explicitly, the Thom isomorphism
\[
\cdot\,t^E(\theta)\colon E^*(X_\ssp)\stackrel{\cong}{\lra}
E^{*+2k}(\Th(\theta))
\]
is given by forming the relative cup product with $t^E(\theta)$, and
is induced by the relative diagonal map $\delta\colon\Th(\theta)\ra
X_\ssp\wedge\Th(\theta)$.

Proposition \ref{thomstr} applies to Definition \ref{defitthomcplx}
whenever the bundles $\theta_i$ are complex. In this case, the
cohomology algebra $E^*(X_i)$ arises from $E^*(X_0)$ by means of
$i$-fold iterated Thom isomorphisms. For example, if $X_0=*$ and
$\dim_\bC\theta_0=k$, then the first iteration identifies $E^*(X_1)$
with $E_*[t]/(t^2)$, where $t$ lies in $E^{2k}(X_1)$; this, of course,
is because $X_1$ is homeomorphic to $S^{2k}$. Further iterations
require the Chern classes of the $\theta_i$.

In subsequent applications, the $\theta_i$ are line bundles. Over $X$,
the isomorphism classes of such bundles form an abelian group with
respect to tensor product, which is isomorphic to $H^2(X)$ under the
Chern class $c_1(-)$. So any $w\in H^2(X)$ gives rise to a complex
line bundle
\begin{equation}\label{deflamb}
\lambda=\lambda(w)\;\;\text{such that $\;c_1(\lambda(w))=w$}\,;
\end{equation}
it is unique up to isomorphism, and
$\lambda(d_1w)=\lambda(w)^{d_1}$ for any integer $d_1$.
The Thom class $t(\lambda(w)^{d_1})$ lies in
$H^2(\Th(\lambda(w)^{d_1}))$, and the second stage of the iteration
begins with a cohomology class $d_2\,t(\lambda(w)^{d_1})$, for some
integer $d_2$. The corresponding line bundle is
$\lambda(w)^{d_1,d_2}\letbe\lambda(t(\lambda(w)^{d_1}))^{d_2}$ over
$\Th(\lambda(w)^{d_1})$. In this language, the $n$th stage identifies 
the bundle
\begin{equation}\label{itc}
\lambda(w)^{d_1,\dots,d_n}\;\letbe\;
\lambda(t(\lambda(w)^{d_1,\dots,d_{n-1}}))^{d_n}\;\;\;\text{over
$\;\Th(\lambda(w)^{d_1,\dots,d_{n-1}})$}\,.
\end{equation}
In other words, $X_n=\Th(\lambda(w)^{d_1,\dots,d_n})$ is an $n$-fold
iterated Thom space over $X_0=X$, for which
$\theta_0=\lambda(w)^{d_1}$ and $\theta_{n-1}=
\lambda(w)^{d_1,\dots,d_n}$.

It is perfectly acceptable to choose $X_0=*$, in which case $w=0$ and
$\lambda(w)$ is the trivial line bundle $\bC$. In this context,
Corollary \ref{divisprop} may be combined with Al Amrani's proof of 
\cite[I.1(c)]{ala:csw} to provide a homeomorphism
\begin{equation}\label{pcasitthom}
\Pc\;\cong\;\Th\big(\lambda(0)^{q_1\pfm
 ,\pfm\dots\pfm,\pfm q_n}\big)
\end{equation}
for any divisive weight vector $\chi$, where $q_j=\chi_j/\chi_{j-1}$
as in Definition \ref{fidivtodiv}.

In general, $E^*(X_n)$ may be computed by iterating Proposition
\ref{thomstr} and exploiting two consequences of the fact that
$\theta$ has dimension~$1$. Firstly, the $E$-theory Thom class
satisfies
\begin{equation}\label{thch}
t^E(\theta)\;=\;c_1^E(\lambda(t(\theta)))
\end{equation}
in $E^2(\Th(\theta))$, which follows directly from the universal
example $\zeta$ over $\bC P^\infty$. Secondly, for any integer $r$,
the equation
\begin{equation}\label{chthk}
c_1^E(\theta^r)\;=\;[r](c_1^E(\theta))
\end{equation}
holds in $E^2(X)$, where $[r]$ denotes the \emph{$r$-series} of the 
formal group law $F_E$ associated to $x^E$ \cite{haz:fga}.
Thus
\[
[r](u)\;\equiv\;r\pfm u+{\textstyle \frac{1}{2}}\,r(r-1)a^Eu^2\mod
(u^3)
\]
in $E^*\llbracket u\rrbracket$, where $a^E\in E_2$ is the coefficient of $u_1u_2$ in
$F_E(u_1,u_2)$.
\begin{thm}\label{itthompspec}
For any divisive $\chi$, the $E_*$-algebra $E^*(\Pc)$ is isomorphic to
\begin{equation}\label{Eofdiv}
E_*[w_n,\,w_{n-1}w_n,\,\dots,\,w_1w_2\cdots w_n]\,\big/\,J^E,
\end{equation}
where $w_hw_{h+1}\cdots w_n$ lies in $E^{2(n-h+1)}(\Pc)$ for any $h\leq i$, 
and $J^E$ denotes the ideal generated by elements of the form 
\[
\big(w_i-[q_i](w_{i-1})\big)\,w_i\cdots w_n
\]
for $1\leq i\leq n$; also $w_0=0$.
\end{thm}
\begin{rem}\label{noexist}
  The elements $w_i$ do not themselves exist in $E^2(\Pc)$ for any
  $i\neq n$, but appear only in monomials divisible by 
  $w_hw_{h+1}\cdots w_n$ for some $h\leq i$. Nevertheless, the description 
  provided by \eqref{Eofdiv} is notationally convenient, and encodes the 
  product structure by repeated application of the relations in $J^E$.
\end{rem}
\begin{proof}[Proof of 4.8]
  Combine Proposition \ref{thomstr} with \eqref{itc},
  \eqref{pcasitthom}, \eqref{thch}, and \eqref{chthk}. 
  The first stage identifies $E^*(\bP(1,1)_\ssp)$ as
  $E_*[w_1]/(w_1^2)$, 
 where $w_1\letbe t^E\?(\lambda(0))$ and $\lambda(0)=\bC$ over $*$.
 The $n$th stage identifies $E^*(\Pc)$ as a free 
 $E^*(\bP(\chi')_\ssp)$-module on the single generator
\[
w_n\;\letbe\; t^E\big(\lambda(0)^{q_1,\dots,\, q_n}\big)\,,
\]
with the relation $w_n^2\;=\;[q_n](w_{n-1})\,w_n$ of \eqref{multreln1}.
\end{proof}
\begin{exa}\label{intco}
  The formal group law associated to integral cohomology is additive,
  and its $r$-series is given by $[r](u)=ru$ in $\bZ\llbracket u\rrbracket$. 
  So for any $p$-primary weight vector $\pi=(1,1,p^{k_2},\dots,p^{k_n})$,
  \eqref{Eofdiv} identifies $H^*(\Pp)$ with
\begin{equation}\label{cohpp}
\bZ[w_n,w_{n-1}w_n,\,\dots,\,w_1w_2\cdots w_n]\,/\,J\,,
\end{equation}
where $J$ is generated by elements of the form
\begin{equation}\label{HQ}
\big(w_i\;-\;p^{k_i-k_{i-1}}w_{i-1}\big)\,w_i\cdots w_n
\end{equation}
for $1\leq i\leq n$, and $w_0=0$. In fact \eqref{cohpp} is
isomorphic to $\bZ[v_1,\dots,v_n]/I(\pi)$ of \eqref{cohringwps}, where
Kawasaki's ideal $I(\pi)$ is generated by the relations
\begin{equation}\label{mpi}
v_1^j=m(\pi)_jv_j,\sts{where}
m(\pi)_j=p^{(j-1)k_n}\!\big/ p^{k_{n-1}+\dots+k_{n-j+1}}\,,
\end{equation}
for $2\leq j\leq n$ \cite[page 243]{kaw:ctp}. The isomorphism arises
from the bijection of generators $w_{n-j+1}\cdots w_n\leftrightarrow
v_j$, by repeated application of \eqref{HQ}; it is multiplicative
because $w_n^j=\prod_{h=1}^{j-1}p^{k_n-k_{n-h}}w_{n-j+1}\cdots w_n$ by
induction on $j$.
\end{exa}
\begin{rems}\label{1andQ}\hfill
\begin{enumerate}
\item[\bf 1.]
Rationally, Theorem \ref{itthompspec} states that $E\bQ^*(\Pc_\ssp)$ is isomorphic to
\[
E\bQ_*[w_n]\,/(w_n^{n+1})\,.
\]
\item[\bf 2.] 
If $\chi=1$, then $\Pc$ reduces to $\CPn$, and Theorem \ref{itthompspec} identifies $w_n$ with $x^E$, and $w_n^j$ with $w_{n-j+1}\cdots w_n$ for 
every $j\geq 1$. 
\end{enumerate}
\end{rems}
These observations illustrate the homotopy equivalences of Remark 
\ref{rattype}. 

\begin{exa}\label{1234c}
The $2$-primary part of $\bP(1,2,3,4)$ is $\bP(1,2,1,4)$, which is
a $2$-fold iterated Thom space over $\bP(1,1)=\bC P^1$. So
$E^*(\bP(1,1)_\ssp)$ is isomorphic to $E_*[w_1]/(w_1^2)$,
where $w_1$ generates $E^2(\bP(1,1))$. Furthermore,
$\lambda(w_1)\cong\zeta$, and $\bP(1,1,2)$ is homeomorphic to
$\Th(\lambda(w_1^2)$; hence $E^*(\bP(1,1,2)_\ssp)$ is isomorphic
to
\[
E_*[w_2,\,w_1w_2]\,\big/\,J^E_1\,,
\]
where $w_2\letbe c_1^E\?(\lambda(w_1)^{2,1})$ and $J_1^E$ is the
ideal generated by
\[
w_1^2w_2\sands w_2^2-2w_1w_2\,.
\]
Similarly, $\bP(1,2,1,4)$ is homeomorphic to
$\Th(\lambda(w_1)^{2,2})$; so $E^*(\bP(1,2,1,4)_\ssp)$ is
isomorphic to
\[
E_*\big[w_3,\,w_2w_3,\,w_1w_2w_3\big]
\,\big/\,J_2^E\,,
\]
where $w_3\letbe c_1^E\?(\lambda(w_1)^{2,2,1})$, and $J_2^E$ is
generated by 
\[
w_1^2w_2w_3,\;\;
\big(w_2^2-2w_1w_2\big)w_3,\;\;\text{and}\;\;w_3^2-2w_2w_3-a^Ew_2^2w_3\,.
\]
The $3$-primary part of $\bP(1,2,3,4)$ is $\bP(1,1,3,1)$, which is a
Thom space over $\bP(1,1,1)=\bC P^2$. So $E^*(\bP(1,1,1)_\ssp)$ is
isomorphic to $E_*[w_2]/(w_2^3)$, where $w_2$ generates
$E^2(\bP(1,1,1))$ and $w_2^2=w_1w_2$. Moreover,
$\lambda(w_2)\cong\zeta$, and $\bP(1,1,1,3)$ is homeomorphic to
$\Th(\lambda(w_1^{1,3})$; so $E^*(\bP(1,1,1,3)_\ssp)$ is
isomorphic to
\[
E^*[w_3,w_2w_3,w_1w_2w_3]\,/J_3^E\,,
\]
where $w_3\letbe c_1^E\?(\lambda(w_1)^{1,3,1})$, and $J_3^E$ is
generated by
\[
w_1^2w_2w_3,\;\;\big(w_2^2-w_1w_2\big)w_3,\;\;
\text{and}\;\;w_3^2-3w_2w_3-3a^Ew_2^2w_3\,.
\]
\end{exa}

The multiplicative formal group law is associated to complex
$K$-theory and the Conner-Floyd orientation. The coefficient ring is
$K_*\cong\bZ[z,z^{-1}]$, and the element $zx^K\in K^0(\bC P^\infty)$
is represented by the virtual Hopf bundle $\zeta-\bC$. The $r$-series
is induced by the tensor power map $\zeta\mapsto\zeta^r$, and is
therefore given by
\begin{equation}\label{Kks}
[r](u)\;=\;z^{-1}\big((1+zu)^r-1\big)
\end{equation}
in $K_*\llbracket u\rrbracket$, for any integer $r$. 
Al Amrani's results of \cite{ala:ckt} may then be recovered.
\begin{exa}\label{Konesq}
Theorem \ref{itthompspec} and \eqref{Kks} combine to show that, for 
any integer $r$, the $K_*$-algebra $K^*(\bP(1,\dots,1,r)_\ssp)$ is 
isomorphic to 
\begin{equation}\label{Kofonesq}
K_*[w_n,\,w_{n-1}w_n,\,\dots,\,w_1w_2\cdots w_n]\,\big/\,J^K,
\end{equation}
where $J^K$ denotes the ideal generated by elements of the form
\[
\big(w_i-w_{i-1}\big)w_i\cdots w_n\;\;\,\text{for}\;\;
1\leq i\leq n-1\,,
\]
and $\big(w_n-z^{-1}\big((1+zw_{n-1})^r-1\big)\big)w_n$. 
The latter is equivalent to 
\[
w_n^2\;=\;\sum_{s=1}^r\binom{r}{s}z^{s-1}w_{n-s}\cdots w_n\,,
\;\text{where}\;\;w_0=0\,.
\]
\end{exa}

%
%
%
%
%
%
%
%
%

\section{Cohomological reassembly}\label{core}

It is now possible to follow the lead of Theorem \ref{reass} by
reassembling the $E_*$-algebra $E^*(\Pc)$ 
from its constituent components $E^*(\bP({}_{p_i}\chi))$.

For any weight vector $\chi$, recall that $Q(\chi)=\{p_1,\dots,p_m\}$ 
denotes the primes occurring in $\chi$. The decomposition of 
Definition \ref{pprim} may then be expressed as 
$\chi={}_{p_i}\chi\alpha(i)$ for each $1\leq i\leq m$, 
where \smash{$Q(\alpha(i))=Q(\chi)\setminus\{p_i\}$}. It is convenient 
to write $\bZ_\chi$ for the subring $\bZ[p_1^{-1},\dots,p_m^{-1}]<\bQ$.
\begin{prop}\label{ecohchifree}
  The $E_*$-algebra $E^*(\Pc_\ssp)$ is a free $E_*$-module, with one 
  generator in each even dimension $\leq 2n$.
\end{prop}
\begin{proof} 
  Consider the insertion map $e(\alpha(j))\colon\bP({}_{p_j}\chi)\to\Pc$ of 
  \eqref{defees}, for some $1\leq j\leq m$. By Proposition \ref{eisfg}, it 
  is the orbit map for the action of the finite group $C_{\alpha(j)}$, whose 
  order is divisible by every $p_i$ such that $i\neq j$. It therefore 
  induces an isomorphism 
\begin{equation}\label{fgiso} 
e(\alpha(j))^*\colon H^*(\Pc_\ssp;\bZ_{\alpha(j)})\;\lra\; 
H^*(\bP({}_{p_j}\chi)_\ssp;\bZ_{\alpha(j)})\,.
\end{equation}
Example \ref{intco} shows that the graded abelian group
$H^*(\bP({}_{p_j}\chi)_\ssp)$ is free, with one generator in each even 
dimension $\leq 2n$. So $H^*(\Pc_\ssp)$ contains at most $p_i$ torsion, for 
$1\leq i\neq j\leq m$. Repeating the argument for every $1\leq j\leq m$ in 
turn proves that $H^*(\Pc_\ssp)$ is torsion free, and therefore has one 
generator in each even dimension $\leq 2n$.

Since $E_*$ is also torsion free and even dimensional, the
Atiyah-Hirzebruch spectral sequence for $E^*(\Pc_\ssp)$ collapses, and
the conclusion follows.
\end{proof}
\begin{cor}\label{corindhom}
For any weight vectors $\chi$ and $\sigma$, the induced homomorphism 
\[
e(\sigma)^*\?\otimes 1\colon E^*(\bP(\sigma\chi)_\ssp)\otimes\bZ_\sigma
\lra E^*(\Pc_\ssp)\otimes\bZ_\sigma
\]
is an isomorphism of algebras over $E_*\?\otimes\bZ_\sigma$.
\end{cor}
\begin{proof}
Proposition \ref{ecohchifree} implies that $e(\sigma)^*$ induces 
an isomorphism of $E_2$-terms of Atiyah-Hirzebruch spectral sequences, 
which collapse. It therefore induces the required isomorphism on their 
limits.
\end{proof}
Proposition \ref{ecohchifree} may, of course, be deduced from Kawasaki's 
calculations; as proven above, it follows from the theory of Thom spaces. 
Isomorphism \eqref{fgiso} confirms that $e(\alpha)\colon\Ppc\to\Pc$ and 
$e(1/\alpha)\colon\Pc\to\Ppc$ are mutually inverse $p$-local homotopy 
equivalences, for any $p$ in $Q(\chi)$.

The next step is to identify a cohomological version of Proposition
\ref{pbs}, by applying $E^*(-)$ to the first diagram \eqref{chirhos}. 
\begin{prop}\label{apos} 
If $\sigma$ and $\sigma'$ are coprime, then the diagram 
\begin{equation}\label{four}
\begin{CD}
E^*(\bP(\sigma\sigma')_\ssp)@<{e(1/\sigma')^*}<<
E^*(\bP(\sigma)_\ssp)\\
@A{e(1/\sigma)^*}AA@AA{e(1/\sigma)^*}A\\
E^*(\bP(\sigma')_\ssp)@<<{e(1/\sigma')^*}<
E^*(\bC P^n_\ssp)\\
\end{CD}
\end{equation}
is a pushout square; in other words, the canonical homomorphism
\[
h\colon E^{*}(\bP(\sigma)_\ssp)\otimes_{E^*(\bC P^n_+)}
E^*(\bP(\sigma')_\ssp)\;\lra\;E^{*}(\bP(\sigma\sigma')_\ssp)
\]
is an isomorphism of $E_*$-algebras.
\end{prop}
\begin{proof}
  Corollary \ref{corindhom} ensures that the horizontal and vertical
  homomorphisms of \eqref{four} induce isomorphisms over
  $E_*\?\otimes\bZ_{\sigma'}$ and $E_*\?\otimes\bZ_\sigma$
  respectively. The same is therefore true of the corresponding
  pushout square. Hence $h$ induces an isomorphism of $E_*$-algebras
  over both $E_*\?\otimes\bZ_{\sigma'}$ and
  $E_*\?\otimes\bZ_\sigma$. But
$\sigma$ and $\sigma'$ are coprime, so $h$ is an isomorphism.
\end{proof}
The cohomological version of Proposition \ref{pos} has a similar proof, with arrows reversed.
\begin{prop}\label{apbs}
If $\sigma$ and $\sigma'$ are coprime, then the diagram 
\begin{equation}\label{apbd}
\begin{CD}
E^*(\bP(\sigma\sigma')_\ssp)@>{e(\sigma)^*}>>
E^*(\bP(\sigma')_\ssp)\\
@V{e(\sigma')^*}VV@VV{e(\sigma')^*}V\\
E^*(\bP(\sigma)_\ssp)@>>{e(\sigma)^*}>
E^*(\bC P^n_\ssp)\\
\end{CD}
\end{equation}
is a pullback square; in other words, the canonical homomorphism 
\[
h\colon E^{*}(\bP(\sigma\sigma')_\ssp)\;\lra\;
E^{*}(\bP(\sigma)_\ssp)\times_{E^*(\bC P^n_\ssp)}E^*(\bP(\sigma')_\ssp)
\]
is an isomorphism of $E_*$-algebras.
\end{prop}
\begin{rem}
 Since $e(\sigma)^*$ and $e(\sigma')^*$ are monic, the limit in Proposition 
 \ref{apbs} may be interpreted as an intersection
\begin{equation}\label{illdesc}
E^*(\bP(\sigma)_\ssp)\cap E^*(\bP(\sigma')_\ssp)\;<\;E^*(\bC P^n_\ssp)\,.
\end{equation}
For $p$ and $p'$-primary weight vectors $\pi$ and $\pi'$, this provides an 
illuminating description of $E^*(\bP(\pi\pi')_\ssp)$ as a subalgebra of 
$E^*(\bC P^n_\ssp)$.
\end{rem}

The cohomological version of Theorem \ref{reass} is now within reach, 
essentially by applying $E^*(-)$ to the geometrical proof.
\begin{thm}\label{asscoho}
For any weight vector $\chi$, the $E_*$-algebra $E^*(\Pc_\ssp)$ is
isomorphic to the colimit of the $\cat{cat}[m]$-diagram
\[
\raisebox{6pt}{$\Ddots$}\hspace{-80pt}
\begin{CD}@.E^*(\bP({}_{p_i}\chi)_\ssp)@.\\
@.@AA{e(1/{}_{p_i}\chi)^*}A@.\\
E^*(\bP({}_{p_1}\chi)\ssp)@<<{e(1/{}_{p_1}\chi)^*}<E^*(\bC P^n_\ssp)
@>>{e(1/{}_{p_m}\chi)^*}>E^*((\bP_m)_\ssp)\,,
\end{CD}
\hspace{-84pt}\raisebox{6pt}{$\ddots$}
\]
and the associated universal homomorphisms 
$E^*(\bP({}_{p_i}\chi)_\ssp)\ra E^*(\Pc_\ssp)$
may be identified with $e({}_{p_i}\chi/\chi)^*$ for every $0\leq i\leq m$.
Similarly, $E^*((\Pc)_\ssp)$ is also isomorphic to the limit of the
$\cat{cat}^{op}[m]$-diagram
\[
\raisebox{6pt}{$\Ddots$}\hspace{-76pt}
\begin{CD}
@.E^*(\bP({}_{p_i}\chi)_\ssp)@.\\
@.@VV{e({}_{p_i}\chi)^*}V@.\\
E^*(\bP({}_{p_1}\chi)_\ssp)@>>e({}_{p_1}\chi)^*>E^*(\bC P^n_\ssp)
@<<e({}_{p_m}\chi)^*<E^*(\bP({}_{p_m}\chi)_\ssp)\,,
\end{CD}
\hspace{-88pt}\raisebox{6pt}{$\ddots$}\\
\]
and the associated universal homomorphisms 
$E^*(\Pc_\ssp)\ra E^*(\bP({}_{p_i}\chi)_\ssp)$ may be identified with 
$e(\chi/{}_{p_i}\chi)^*$ for every $0\leq i\leq m$.
\end{thm}
\begin{proof}
  Proceed by induction on $m$, as in the proof of Theorem \ref{reass}.
  The inductive steps appeal to Propositions \ref{apos} and \ref{apbs}
  respectively.
\end{proof}
\begin{rem}
Theorem \ref{asscoho} shows that $E^*(-)$ converts the geometric
limits and colimits of Theorem \ref{reass} into the
corresponding algebraic colimits and limits. Although the geometric
pullbacks are 
not of fibrations, the induced algebraic pushouts are those of a 
collapsed Eilenberg-Moore spectral sequence.
\end{rem}

The pullback and pushout descriptions of Theorem \ref{asscoho} yield
isomorphisms
\begin{equation}\label{ecopb}
E^*(\Pc_\ssp)\;\stackrel{\cong}{\lra}\;
E^*(\bP({}_{p_1}\chi)_\ssp)\otimes_{E^*(\bC P^n_\ssp)}\dots
\otimes_{E^*(\bC P^n_\ssp)}E^*(\bP({}_{p_m}\chi)_\ssp)
\end{equation}
and
\begin{equation}\label{ecopo}
E^*(\bP({}_{p_1}\chi)_\ssp)\times_{E^*(\bC P^n_\ssp)}\dots
\times_{E^*(\bC P^n_\ssp)}E^*(\bP({}_{p_m}\chi)_\ssp)
\;\stackrel{\cong}{\lra}\;E^*(\Pc_\ssp)
\end{equation}
respectively; by analogy with \eqref{illdesc}, the latter may be 
rewritten as
\begin{equation}\label{ecoint}
E^*(\bP({}_{p_1}\chi)_\ssp)\cap\dots\cap E^*(\bP({}_{p_m}\chi)_\ssp)\;<\; 
E^*(\bC P^n_\ssp)\,.
\end{equation}
Each of these $E_*$-algebras has one $E_*$-generator in each even 
dimension $\leq 2n$, and leads directly back to Kawasaki's original calculations.
\begin{thm}\label{kawa}
In the case of integral cohomology, the isomorphisms \eqref{ecopb},
\eqref{ecopo}, and \eqref{ecoint} identify $H^*(\Pc_\ssp)$ with
$\bZ[v(1),\dots,v(n)]/I(\chi)$, as in \eqref{cohringwps}.
\end{thm}
\begin{proof}
  For each prime $p_i\in Q(\chi)$, let $v_i\in H^{2j}(\Pc_\ssp)$
  denote the image of Kawasaki's generator $1\otimes\dots\otimes
  v_i(j)\otimes\dots\otimes 1$ under the isomorphism \eqref{ecopb},
  where $v_i(j)\in H^{2j}(\bP({}_{p_i}\chi))$. Example \ref{intco} shows 
  that $m_i(j)v_i(j)=v_i(1)^j$ lies in the image of $H^{2j}(\bC P^n)$ for
  every $i$, where
\[
m_i(j)\;=\;p_i^{(j-1)k_i(n)}/p_i^{k_i(n-1)+\dots+k_i(n-j+1)}
\]
as in \eqref{mpi}; thus $m_1(j)v_1=m_i(j)v_i$ for every $i$. The
numbers $M_i\letbe\prod_{h\neq i}m_h(j)$ are coprime, and satisfy
$M_iv_h=M_hv_i$ for every $h\neq i$. So there exist non-zero integers
$A_i$ such that $\sum_iA_iM_i=1$. Now define the element $v(j)$ to be
$\sum_iA_iv_i$ in $H^{2j}(\Pc)$. For every $h$, it follows that
\[
M_hv(j)\;=\;\sum_iA_iM_hv_i\;=\;\big(\sum_iA_iM_i\big)v_h\;=\;v_h,
\]
and hence that $H^{2j}(\Pc)$ is free abelian, on generator $v(j)$. By
construction, $v(1)^j=\prod_im_i(j)v(j)$ for every $1\leq j\leq n$, as
required by \eqref{cohringwps}.

Alternatively, the isomorphism \eqref{ecoint} identifies 
$v_i(j)\in H^{2j}(\bP({}_{p_i}\chi)$ with the element
$p_i^{k_i(n)+\dots+k_i(n-j+1)}x^j\in H^{2j}(\bC P^n)$. The intersection of 
the cyclic groups so generated is therefore infinite cyclic on 
\[
v(j)\;=\;\prod_ip_i^{k_i(n)+\dots+k_i(n-j+1)}x^j\,,
\]
and the relation $v(1)^j=\prod_im_i(j)v(j)$ follows again. 
\end{proof}

\begin{exa}\label{345c}
By \eqref{ecoint}, $K^*(\bP(3,4,5)_\ssp)$ may be identified with
\begin{equation}\label{k345f}
K^*(\bP(1,4,1)_\ssp)\cap K^*(\bP(3,1,1)_\ssp)\cap K^*(\bP(1,1,5)_\ssp)\;<\; 
K^*(\bC P^2_\ssp)\,.
\end{equation}
Example \ref{Konesq} with $n=2$ shows that 
$K^*(\bP(1,1,r)_\ssp)$ is isomorphic to 
\begin{equation}\label{K11q}
K_*[w_2,\,w_1w_2]\,\big/\,
\big(w_1^2w_2,\,w_2^2-rw_1w_2\big)\,,
\end{equation}
and that \eqref{ecoint} identifies $w_2$ with $[r](x^K)$ in 
$K^2(\bC P^2)$ and $w_1w_2$ with \smash{$r(x^K)^2$} in 
$K^4(\bC P^2)$. Substituting $r=3$, $4$, and $5$ in turn confirms 
that $K^*(\bP(3,4,5))$ has
$y_1=60x^K+90z(x^K)^2$ and $y_2=60(x^K)^2$ as $K_*$-generators, 
in dimensions $2$ and $4$ respectively. So there is an isomorphism
\begin{equation}\label{K345e}
K^*(\bP(3,4,5)_\ssp)\;\cong\;
K_*[y_1,y_2]\,\big/\,(y_1^2-60y_2)\,.
\end{equation}
\end{exa}
Setting $z=1$ in \eqref{K345e} provides an example of Al Amrani's abstract 
isomorphism \cite[Corollary 3.2]{ala:cbc} between $K^*(\Pc)$ and 
$H^*(\Pc)$ for certain $\chi$.

Theorem \ref{kawa} expresses Kawasaki's discovery (which he did not
make explicit) that the ring \smash{$H^*(\Pc_\ssp)$} depends only on
the unordered coordinates of the vectors ${}_p\chi$, as $p$ ranges
over $Q(\chi)$. The same holds for the $E_*$-algebra $E^*(\Pc_\ssp)$.

\begin{rem}\label{htpytype2}
These facts also follow from \cite{bafrnora:cwp}. As explained in Remark 
\ref{htpytype}, $E^*(\Pc)$ may always be described in terms of iterated 
Thom isomorphisms as $E^*(\bP(\chi^*))$. The advantage presented by 
Theorem \ref{asscoho} is that the $p$-primary parts $E^*(\bP({}_p\chi))$ are 
each computed using Theorem \ref{itthompspec}; since the computations 
involve the $E$-theory $p^k$-series (rather than the $r$-series for 
composite $r$), the technical machinery of Brown-Peterson cohomology theory 
\cite{rav:ccs} may then be brought to bear. 
\end{rem}

%
%
%
%
%
%
%
%
%

\section{Homology and fundamental classes}\label{hofucl}

Since \daaja's original work \cite{daja:cpc}, toric topology has
tended to focus on cohomological calculations to the detriment of
their homological counterparts. For weighted projective spaces,
however, the complex bordism coalgebras $\osu(\Pc)$ are of particular
interest, and this section is devoted to understanding $E_*(\Pc)$ for
any complex oriented homology theory $E_*(-)$. 

For $\CPn$, the complex orientation reveals itself as an isomorphism
\begin{equation}\label{ehpn}
E_*(\bC P^n_\ssp)\;\stackrel{\cong}{\llra}\;E_*\br{b_0,b_1,\dots,b_n}
\end{equation}
of free $E_*$-coalgebras, where $b_j$ has dimension $2j$ and
supports the coproduct
\[
\delta(b_j)\;=\;\sum_{i=0}^jb_i\otimes b_{j-i}
\]
in $E_*(\bC P^n_\ssp)\otimes E_*(\bC P^n_\ssp)$; the $b_j$ form the
dual $E_*$-basis to the powers \smash{$(x^E)^j$} for $0\leq j\leq n$,
and $b_0$ is the counit $1$. 

For notational clarity, two conventions are adopted throughout the 
remainder of this section. Firstly, $b_j$ is expanded to $b^E_j$ 
whenever the homology theory needs emphasising; and secondly, following 
Chapter \ref{intro}, the universal complex orientation is usually denoted 
by $u$ in $\varOmega^2(\bC P^n)$.

It is important to clarify the relationship between the $b_j$ and the 
Poincar\'e duality isomorphism
\begin{equation}\label{epdi}
{}\cap\sigma\colon E^i(\bC P^n_\ssp)\;\stackrel{\cong}{\llra}\;
E_{2n-i}(\bC P^n_\ssp)\,,
\end{equation}
defined by cap product with a fundamental class 
$\sigma\in E_{2n}(\bC P^n_\ssp)$. This is best done in the context of the
universal example, and has an interesting history.

During the early days of the theory, it was usual to identify
$\osu(\bC P^n_\ssp)$ with the free $\osu$-module on generators
$cp_j$, represented by the inclusions $\bC P^j\to\CPn$ for $0\leq
j\leq n$. From this viewpoint, $cp_n$ is the bordism class of the
identity map $1_{\bC P^n}$, and the most natural choice of
fundamental class $\sigma$. In particular, iteration of the
formula
\begin{equation}\label{xcpn}
u\pfm\cap cp_n\;=\;cp_{n-1}
\end{equation}
in $\osu(\bC P^n_\ssp)$ shows that $cp_j$ is the Poincar\'e dual of
$u^{n-j}$ for every $0\leq j\leq n$. On the other hand, the
$cp_j$ are certainly not $\Hom$ dual to the $u^j$, but may
be expanded by
\begin{equation}\label{cpbj}
cp_j\;=\;[\bC P^j]1+[\bC P^{j-1}]b_1+\dots+[\bC P^1]b_{j-1}+b_j
\end{equation} in terms of the basis \eqref{ehpn}.

Formula \eqref{cpbj} is originally due to Novikov, and is
an immediate consequence of \eqref{xcpn}. It emphasises the fact
that $b_j$ lies in the \emph{reduced} group 
\smash{$\varOmega_{2j}^U(\bC P^n)$} for every $1\leq j\leq n$, 
whereas $cp_j$ has obvious non-trivial augmentation. Nevertheless, 
$b_n$ may be deployed equally well as
a fundamental class, and determines an alternative Poincar\'e
duality isomorphism. Since
\[
u\pfm\cap b^U_n\;=\;b^U_{n-1} 
\]
holds in $\osu(\CPn)$ by definition, $b_j$ is the alternative
Poincar\'e dual of $u^{n-j}$ by analogy with \eqref{xcpn}. Under the
Thom orientation $\varOmega^U_*(-)\to H_*(-)$, both $b_n$ and $cp_n$
map to the canonical fundamental class $b^H_n$ in $H_{2n}(\CPn)$. They
therefore induce the same Poincar\'e duality isomorphism in integral
homology and cohomology.

The problem arises of identifying geometrical representatives
$B_j\ra\CPn$ for the $b_j$, bearing in mind that the stably
complex manifolds $B_j$ must bound when $j\geq 1$, because the $b_j$
are reduced. This was solved in \cite{ray:cbt}, where the $B_j$ are
constructed as iterated sphere bundles whose stably complex
structures extend over the associated disc bundles. Subsequently,
the $B_j$ were identified as \emph{Bott towers} \cite{grka:btc}, and
therefore as non-singular toric varieties with canonical complex
structures; the stabilisations of these structures do \emph{not}
bound, having non-trivial Chern numbers. The language of iterated
sphere bundles is documented in \cite[Sections 2 \& 3]{cira:hdr}, and 
used extensively below.

Each Bott tower is determined by a list $(r_1,\dots,r_n)$ of integral
$j$--vectors $r_j$. Given any divisive weight vector $\chi$, let
$q_j=\chi_j/\chi_{j-1}$ as in Definition \ref{fidivtodiv}, and define
the Bott tower $(B_j(\chi):0\leq j\leq n)$ by choosing
\[
r_j\;=\;(0,\dots,0,q_j)\,.
\]
Then $B_0(\chi)=*$, and $B_j(\chi)$ is a $2j$-dimensional stably
complex manifold equipped with canonical complex line bundles
$\gamma_i=\gamma_i(\chi)$, for $0\leq i\leq j\leq n$. It is defined
inductively as the total space $S(\delta_j(\chi))$ of the
$2$--sphere bundle of
\begin{equation}\label{defdel}
\delta_j(\chi)\;\letbe\;\bR\oplus\gamma_{j-1}^{q_j}
\end{equation}
over $B_{j-1}(\chi)$, where $\bR$ denotes the trivial real line bundle.
The unit $1\in\bR$ determines a section $i_{j-1}$ for
$\delta_j(\chi)$, which features in the cofibre sequence
\begin{equation}\label{cofib}
B_{j-1}(\chi)\stackrel{i_j}{\llra}
B_j(\chi)\stackrel{l_j}{\llra}\Th(\gamma_{j-1}^{q_j})\,.
\end{equation}
In terms of the complex orientation $x^E$, the corresponding Thom
class \smash{$t^E(\gamma_{j-1}^{q_j})$} generates
\smash{$E^2(\Th(\gamma_{j-1}^{q_j}))$} and pulls back to the generator
$v_j=v_j^E$ of $E^2(B_j(\chi))$, as described in \cite[Chapter
3]{cira:hdr}. The inductive description is completed by appealing to
\eqref{deflamb}, and letting $\gamma_j$ be the complex line bundle
$\lambda(v^H_j)$.

The stably complex structure on $B_j(\chi)$ is induced from the
defining $S^2$-bundle \eqref{defdel}, and extends over the $3$-disc
bundle $D(\delta_j(\chi))$ by Szczarba \cite{szc:tbf}. It is specified 
by a canonical isomorphism
\begin{equation}\label{scsbpi}
c\colon\tau(B_j(\chi))\oplus\bR\stackrel{\cong}{\llra}
\gamma_0^{q_1}\!\oplus\cdots
\oplus\gamma_{j-1}^{q_j}\!\oplus\bR
\end{equation}
of $SO(2j+1)$-bundles, where $\gamma_0=\bC$.

Recall that \eqref{pcasitthom} expresses $\Pc$ as an iterated Thom
space over $*$. The sequence $X_0$, $X_1$, \dots, $X_n$ of Thom
spaces may be written as
\begin{equation}\label{thomseq1}
*,\;\;\Th(\lambda(0)^{q_1}),\;\;\Th(\lambda(0)^{q_1,\pfm q_2}),\;\;\dots,\;\;
\Th(\lambda(0)^{q_1,\pfm q_2,\pfm\dots,\pfm q_n})\,,
\end{equation}
or equivalently as
\begin{equation}\label{thomseq2}
\bP(1),\;\;\bP(1,\chi_1),\;\;\bP(1,\chi_1,\chi_2),\;\;\dots,\;\;
\bP(1,\chi_1,\chi_2,\dots,\chi_n)\,.
\end{equation}
For any $1\leq j\leq n$,
the Thom class $t^E(\lambda(0)^{q_1,\pfm\dots,\pfm q_j})$ that arises
from \eqref{thomseq1} coincides with the generator $w_j$ in
$E^2(\bP(1,\chi_1,\dots,\chi_j))$
that arises from \eqref{thomseq2}; it is convenient to 
denote them both by $t_j$ in $E^2(X_j)$.

\begin{lem}\label{fjn}
For every $1\leq j\leq n$, there exists a map 
$f_{j,n}\colon B_j(\chi)\to X_n$ such that 
$f_{j,n}^*(t_n)=v_j$ in $E^2(B_j(\chi))$.
\end{lem}
\begin{proof}
Proceed by induction on $j$, with base case $j=1$. 

For any $n\geq 1$, the map $f_{1,n}\colon B_1(\chi)\to X_n$ is
necessarily the inclusion of the fibre
$S^2\subset\Th(\lambda(0)^{q_1,\dots,q_n})$, and coincides with the
map $\bP(\chi_{n-1},\chi_n)\to\Pc$ induced on the final two
homogeneous coordinates.

Assume that $f_{j-1,n}$ exists with the required properties, and choose 
$2\leq j\leq n$. Thus $1\leq j-1\leq n-1$, and it follows that 
$f_{j-1,n-1}\colon B_{j-1}(\chi)\to X_{n-1}$ satisfies 
$f_{j-1,n-1}^*(t_{n-1})=v_{j-1}$ in $E^2(B_{j-1}(\chi))$. Now define 
$f_{j,n}$ as the composition
\begin{equation}\label{fjndef}
B_j(\pi)\stackrel{l_j}{\lllllra}\Th(\gamma_{j-1}^{q_j})
\stackrel{f'_{\?j-1,n-1}}{\lllllra}X_n
\end{equation}
for $2\leq j\leq n$, where $\gamma_{j-1}=\lambda(v_{j-1}^H)$ and
$f'_{j-1,n-1}$ denotes $\Th(f_{j-1,n-1})$. Thus
$(f_{j-1,n-1}')^*(t_n)=t^E(\gamma_{j-1}^{q_j})$ in
$E^2(\Th(\gamma_{j-1}^{q_j}))$, and applying $l_j^*$ yields the
required equation.
\end{proof}

Every complex orientation induces natural transformations 
$\osu(-)\ra E_*(-)$ and $\ous(-)\ra E^*(-)$, both of which are written  
$x^E_*$. They reduce to the identity in the universal case.
\begin{defn}\label{defbujpi}
For any $0\leq j\leq n$, the bordism class $b_j(\chi)$ is represented by 
the map $f_{j,n}$ of Lemma \ref{fjn}, and lies in 
$\varOmega^U_{2j}(\Pc_\ssp)$; its image $x^E_*(b_j(\chi))$ is also denoted 
by $b_j(\chi)$ (or $b_j^E(\chi)$ to avoid ambiguity), 
and lies in $E_{2j}(\Pc_\ssp)$. 
\end{defn}
The $B_j(\chi)$ bound as stably complex manifolds for every $j>0$, so
the corresponding $b_j(\chi)$ actually belong to the reduced groups 
$E_{2j}(\Pc)$.

According to Theorem \ref{itthompspec}, the elements 
\begin{equation}\label{ouspibasis}
\{w_{i+1}\cdots w_n:0\leq i\leq n\}
\end{equation}
form an $\osu$-basis for $\ous(\Pc_\ssp)$, where the case $i=n$ is  
interpreted as $1$. 

\begin{thm}\label{osuppi}
For any divisive $\chi$, the elements 
$\{b_j(\chi):0\leq j\leq n\}$ form a basis for the 
$\osu$-coalgebra $\osu(\Pc_\ssp)$; this basis is dual to \eqref{ouspibasis}.
\end{thm}
\begin{proof}
Proceed by induction on $n$, noting that the result is trivial for $n=0$ and 
$\chi=(1)$. The inductive assumption is that 
$\{b_j(\chi'):0\leq j\leq n-1\}$ and 
\[
\{w_{i+1}\cdots w_{n-1}:0\leq i\leq n-1\}
\]
are dual $\osu$-bases for $\osu(\bP(\chi')_\ssp)$ and 
$\ous(\bP(\chi')_\ssp)$ respectively. In other words, the Kronecker product 
$\bbr{w_{i+1}\cdots w_{n-1},\,b_j(\chi')}$ evaluates to 
$\delta_{n-i-1,j}$ in $\varOmega^U_{2(j-i)}$, for every $0\leq i,j\leq n-1$.

Now consider \eqref{fjndef}, and the Thom isomorphisms 
\[
\varPhi^*\colon\ous((X_{n-1})_\ssp)\lra\varOmega^{*+2}_U(X_n) 
\sands
\varPhi_*\colon\varOmega_{*+2}^U(X_n)\lra\osu((X_{n-1})_\ssp)
\]
given by the Thom class $w^U\?(n)\in\varOmega^2_U(X_n)$. Then
$\varPhi^*$ satisfies
\begin{equation}\label{cobthomiso}
\varPhi^*(w_{i+1}\cdots w_{n-1})\;=\;w_{i+1}\cdots w_n
\end{equation} 
for $0\leq i\leq n-1$, and 
Lemma \ref{fjn} shows that $\varPhi_*$ satisfies
\begin{equation}\label{borthomiso}
\varPhi_*(b_k(\chi))\;=\;b_{k-1}(\chi')
\end{equation}
for $1\leq k\leq n$. To check that the stably complex structures
behave as required by \eqref{borthomiso}, appeal must be made to 
\eqref{scsbpi}.

By \eqref{cobthomiso} and \eqref{borthomiso}, the Kronecker product 
$\bbr{w_{i+1}\cdots w_n,\,b_k(\chi)}$ may be rewritten as
\[
\begin{split}
\bbr{\varPhi^*(w_{i+1}\cdots w_{n-1}),\,b_k(\chi)}
\;&=\;
\bbr{w_{i+1}\cdots w_{n-1},\,\varPhi_*(b_k(\chi))}\\
&=\;\bbr{w_{i+1}\cdots w_{n-1},\,b_{k-1}(\chi')}
\end{split}
\]
for $0\leq i\leq n-1$ and $1\leq k\leq n$. This evaluates to 
$\delta_{n-i-1,k-1}=\delta_{n-i,k}$ by induction, so it remains only to 
confirm the cases $i=n$ and $k=0$. They involve 
$1\in\varOmega^0_U((X_n)_\ssp)$ and $1\in\varOmega_0^U((X_n)_\ssp)$, and 
follow immediately.
\end{proof}

\begin{cor}\label{eppi}
  Given any complex oriented homology theory $E_*(-)$, the elements
  $\{b_j(\chi):0\leq j\leq n\}$ form an $E_*$-basis for
  $E_*(\Pc_\ssp)$; it is dual to the basis
  $\{w_{i+1}\cdots w_n:0\leq i\leq n\}$ for $E^*(\Pc_\ssp)$
  given by Theorem \ref{itthompspec}.
\end{cor}
It follows from Corollary \ref{eppi} that the coalgebra structure on 
$E_*(\Pc_\ssp)$ is determined by expressions of the form 
\begin{equation}\label{delbe}
\delta(b_j(\chi))\;=\;
\sum_{0\leq k+l\leq j}e_{j,k,l}\,b_k(\chi)\otimes b_l(\chi)
\end{equation}
in $E_*(\Pc_\ssp)\otimes E_*(\Pc_\ssp)$, for $0\leq j\leq n$. The 
coefficient $e_{j,k,l}$ lies in $E_{2(j-k-l)}$, and is given by the 
coefficient of $w_{n-j+1}\cdots w_n$ in the product
\[
w_{n-k+1}\cdots w_n\cdot w_{n-l+1}\cdots w_n.
\]

\begin{exa}\label{E111q}
If $\chi=(1,1,1,r)$, then Example \ref{1234c} shows that 
$E_*(\Pc_\ssp)$ is freely generated over $E_*$ by elements 
$\{b_j\letbe b_j(\chi):0\leq k\leq 2\}$, where $b_0=1$. Applying 
Theorem \ref{itthompspec} to \eqref{delbe} yields $\delta(1)=1\otimes 1$,
together with
\[
\begin{split}
\delta(&b_1)\;=\;b_1\otimes 1+1\otimes b_1,\quad
\delta(b_2)\;=\;b_2\otimes 1+rb_1\otimes b_1+1\otimes b_2,\\
\sands &\delta(b_3)\;=\;b_3\otimes 1+rb_2\otimes b_1+
\binom{r}{2}a^Eb_1\otimes b_1+rb_1\otimes b_2+1\otimes b_3 
\end{split}
\]
in $E_*(\Pc_\ssp)\otimes E_*(\Pc_\ssp)$.
\end{exa}

In the case of integral homology, the top dimensional group
$H_{2n}(\Pc)$ is isomorphic to $\bZ$, and cap product with the
generator $b^H_n(\chi)$ defines Poincar\'e duality over $\bQ$; this
does not, of course, extend to $\bZ$, symptomising the existence of
singularities.  Nevertheless, $b^H_n(\chi)$ may still be thought of as
a fundamental class, and is the image of the universal $b_n(\chi)$
under $x^H_*$. In this sense, the representing map $j_{n,n}\colon
B_n(\chi)\ra\Pc$ may be interpreted as a desingularisation of $\Pc$;
it is closely related to the associated toric desingularisation
\cite[\S 2.6]{ful:itv}, as we shall explain in a future note.  


%
%
%
%
%
%
%
%
%

\section{Homological reassembly}\label{hore}

It remains to consider the assembly problem for $E_*(\Pc)$, by
dualising the results of Section \ref{core}. This approach is not
strictly necessary, as explained in Remark \ref{htpytype2}; for
$E_*(\Pc_\ssp)$ is isomorphic as $E_*$-coalgebras to
$E_*(\bP(\chi^*))$, and the latter may be described in terms of
iterated Thom isomorphisms.  Nevertheless, the homological advantages
of proceeding prime by prime are as valid as for cohomology.

\begin{prop}\label{ehohchifree}
For any weight vector $\chi$, the $E_*$-coalgebra $E_*(\Pc_\ssp)$ is  a 
free $E_*$-module, with one generator in each even dimension $\leq 2n$.
\end{prop}
\begin{proof}
Because $E_*$ is even dimensional and torsion free, the result follows 
directly from applying $\Hom_{E_*}(-,E_*)$ to Proposition  
\ref{ecohchifree}.
\end{proof}

A more explicit description is obtained by dualising the free 
$E_*$-modules that appear in Theorem \ref{asscoho}.
\begin{thm}\label{assho} 
For any weight vector $\chi$, the $E_*$-coalgebra $E_*(\Pc_\ssp)$ is 
isomorphic to the limit of the $\cat{cat}^{op}[m]$-diagram
\[
\raisebox{6pt}{$\Ddots$}\hspace{-84pt}
\begin{CD}
@.E_*(\bP({}_{p_i}\chi)_\ssp)@.\\
@.@VV{e(1/{}_{p_i}\chi)_*}V@.\\
E_*(\bP({}_{p_1}\chi)_\ssp)@>>{e(1/{}_{p_1}\chi)_*}>E_*(\bC P^n_\ssp)
@<<{e(1/{}_{p_m}\chi)_*}<E_*(\bP({}_{p_m}\chi)_\ssp)\,;
\end{CD}
\hspace{-88pt}\raisebox{6pt}{$\ddots$}
\]
the corresponding universal map $E_*(\Pc_\ssp)\ra E_*(\bP({}_{p_i}\chi)_\ssp)$ 
may be identified with $e({}_{p_i}\chi/\chi)_*$ for every $0\leq i\leq m$. 
Similarly, $E_*(\Pc_\ssp)$ is also isomorphic to the colimit of the
$\cat{cat}[m]$-diagram
\[
\raisebox{6pt}{$\Ddots$}\hspace{-80pt}
\begin{CD}
@.E_*(\bP({}_{p_i}\chi)_\ssp)@.\\
@.@AA{e({}_{p_i}\chi)_*}A@.\\
E_*(\bP({}_{p_1}\chi)_\ssp)@<<{e({}_{p_1}\chi)_*}<
E_*(\bC P^n_\ssp)@>>{e({}_{p_m}\chi)_*}>E_*(\bP({}_{p_m}\chi)_\ssp)\,;
\end{CD}
\hspace{-84pt}\raisebox{6pt}{$\ddots$}
\]
the corresponding universal map $E_*(\bP({}_{p_i}\chi)_\ssp)\ra E_*(\Pc_\ssp)$ 
may be identified with $e(\chi/{}_{p_i}\chi)_*$ for every $0\leq i\leq m$.
\end{thm}

The limit described by Theorem \ref{assho} is actually the iterated 
pullback
\begin{equation*}\label{ehoitpb}
E_*(\bP({}_{p_1}\chi)_\ssp)\times_{E_*(\bC P^n_\ssp)}\dots
\times_{E_*(\bC P^n_\ssp)}E_*(\bP({}_{p_m}\chi)_\ssp)
\end{equation*}
of $E_*$-coalgebras, and the colimit is the iterated pushout
\begin{equation*}\label{ehoitpo}
E_*(\bP({}_{p_1})_\ssp)\otimes_{E_*(\bC P^n_\ssp)}\dots
\otimes_{E_*(\bC P^n_\ssp)}E_*(\bP({}_{p_m})_\ssp)\,.
\end{equation*}
By analogy with \eqref{ecoint}, the former may be rewritten as 
\begin{equation}\label{ehoint}
E_*(\bP({}_{p_1}\chi)_\ssp)\cap\dots\cap E_*(\bP({}_{p_m}\chi)_\ssp)\;<\; 
E_*(\bC P^n_\ssp)\,.
\end{equation}

\begin{exa}\label{345d}
Expression \eqref{ehoint} identifies $K_*(\bP(3,4,5)_\ssp)$ with
\begin{equation}\label{k345}
K_*(\bP(1,4,1)_\ssp)\cap K_*(\bP(3,1,1)_\ssp)\cap K_*(\bP(1,1,5)_\ssp)
\;<\;K_*(\bC P^2_\ssp)\;.
\end{equation}
Applying Corollary \ref{eppi} with $\pi=(1,1,r)$ 
shows that $K_*(\bP(\pi)_\ssp)$ is isomorphic to 
\begin{equation}\label{Kho11q}
K_*\big\langle 1,b_1(\pi),b_2(\pi)\big\rangle\,,
\end{equation}
and also that \eqref{ehoint} identifies $b_1(\pi)$ with 
$b_1$ in $K_2(\bC P^2)$ and $b_2(\pi)$ with $rb_2$ in $ K_4(\bC P^2)$.
These identifications are compatible with the diagonals 
\begin{multline*}
\delta(b_1(\pi))\;=\;b_1(\pi)\otimes 1+1\otimes b_1(\pi)\sands\\
\delta(b_2(\pi))\;=
\;b_2(\pi)\otimes 1+rb_1(\pi)\otimes b_1(\pi)+1\otimes b_2(\pi)\,.
\end{multline*}
Substituting $r=3$, $4$, and $5$ and applying \eqref{ehoint} confirms 
that $K_*(\bP(3,4,5)_\ssp)$ is isomorphic to the subcoalgebra 
\begin{equation}\label{345subco}
K_*\big\langle 1,b_1,60b_2\big\rangle\;<\;K_*(\bC P^2_\ssp)\,,
\end{equation}
as predicted by Proposition \ref{ehohchifree}. The resulting coalgebra is 
$K_*$-dual to the $K_*$-algebra description of $K^*(\bP(3,4,5)_\ssp)$ 
given by Example \ref{345c}.
\end{exa}
Setting $z=1$ in \eqref{345subco} yields an abstract isomorphism of
coalgebras between $K_*(\Pc)$ and $H_*(\Pc)$, and dualises Al Amrani's 
algebra isomorphism of \cite{ala:cbc}.

%
%
%
%
%
%
%
%
%


\end{document}